\documentclass{ifacconf}

\usepackage{natbib}
\usepackage{graphicx}      
\usepackage{amsmath}
\usepackage{amssymb}
\usepackage{enumitem}
\usepackage{mathtools}
\usepackage{xcolor}
\usepackage{tikz,pgfplots}
\usepackage[american]{circuitikz}
\usepackage{anyfontsize}
\pgfplotsset{compat=1.18} 

\allowdisplaybreaks

\newtheorem{objective}{\textbf{Objective}}
\newtheorem{theorem}{Theorem}
\newtheorem{lemma}{Lemma}
\newtheorem{corollary}{Corollary}
\newtheorem{assumption}{Assumption}
\newtheorem{remark}{Remark}
\newtheorem{proposition}{Proposition}
\newtheorem{example}{Example}

\newcommand{\ee}{\mathrm{e}}
\newcommand{\R}{\mathbb{R}}
\newcommand{\N}{\mathbb{N}}

\DeclareMathOperator*{\argmin}{argmin}

\newcommand{\GG}[1]{\textcolor{blue}{#1}}
\newcommand{\RK}[1]{{\color{blue}{#1}}}
\newcommand{\AR}[1]{\textcolor{blue}{#1}}

\begin{document}
\begin{frontmatter}

\title{{Observer-Based Sampled-Data Stabilisation of Switched Systems with Lipschitz Nonlinearities and Dwell-Time}\thanksref{footnoteinfo}}
\thanks[footnoteinfo]{R. Katz is supported by the Alon Fellowship from the Council of Higher Education of Israel. G. Giordano acknowledges support from the European Union through the ERC INSPIRE grant (project number 101076926); views and opinions expressed are however those of the authors only and do not necessarily reflect those of the European Union or the European Research Council; neither the European Union nor the granting authority can be held responsible for them.
}

\author[First]{Rami Katz} 
\author[Second]{Antonio Russo} 
\author[Third]{Gian Paolo Incremona}
\author[Third]{Patrizio Colaneri}
\author[Fourth]{Giulia Giordano}

\address[First]{School of Electrical and Computer Engineering, Tel Aviv University, Israel (e-mail: ramkatsee@tauex.tau.ac.il).}
\address[Second]{Dipartimento di Ingegneria Gestionale, dell’Informazione e della Produzione, University of Bergamo, Italy (e-mail: antonio.russo@unibg.it)}
\address[Third]{Dipartimento di Elettronica, Informazione e Bioingegneria, 
   Politecnico di Milano, Italy (e-mails: \{gianpaolo.incremona,patrizio.colaneri\}@polimi.it)}
\address[Fourth]{Dipartimento di Ingegneria Industriale, University of Trento, Italy (e-mail: giulia.giordano@unitn.it).}

\begin{abstract}                
We investigate the stabilisation of nominally linear-affine switched systems with uncertain Lipschitz nonlinearities under dwell-time constraints, using a sampled-data switching law based on a state observer. 
We design the switching law based on Lyapunov-Metzler inequalities, accounting for the sampled-data output measurements, and we derive time-dependent LMI conditions for global asymptotic stability (or, in the presence of switching affine terms, ultimate boundedness) of the resulting closed-loop system. We obtain an estimate of the average quadratic cost and a bound on its maximum deviation from the actual cost.
Moreover, we discuss the feasibility of the derived LMIs. Specifically, we show how the observer gains can be incorporated into the matrix inequalities, provide equivalent reduced-order LMI conditions, and prove that the time dependence of the LMIs can be removed by discretising on a finite grid. Numerical examples, including practical applications to real-world engineering scenarios in power systems, illustrate our theoretical results and compare them with an existing approach for output-feedback stabilisation of switched systems, subject to sampled-data measurements.
\end{abstract}


\end{frontmatter}

\section{Introduction}
Switched systems are composed of multiple modes, each governed by a distinct dynamical system, while a switching law dictates the transitions between modes \citep{switching:Liberzon:2003}. The switched system does not necessarily retain the stability characteristics of the individual modes; hence, it is crucial to analyse the stability of the switched system under various switching conditions \citep{switched:LiberzonMorse:1999} and to design stabilising switching strategies \citep{LM:GeromelColaneri:2006}.
Classic design approaches rely either on arbitrarily fast switching among modes or on a dwell time, i.e., a minimum amount of time between consecutive switching instants \citep{LM:GeromelColaneri:2006}. Particular attention has been devoted to switched systems that are either linear or affine, in continuous time \citep{ctswitched:Rantzer:1998,ctswitched:Hespanha:2004,quadratic:bolzern:2004,LM:GeromelColaneri:2006,time-triggered:albea:2021, affine:Russo:2024} and in discrete time \citep{dtswitched:Daafouz:2001,discrete:Geromel:2006,LM:Egidio:2019,outputfb:Egidio:2021}.
While the literature mostly considers full state feedback, a body of work focuses on the design of stabilising output feedback laws \citep{outpput:Faron:1996,output:DeCarlo:2000,outputfeedback:Deaecto:2016,outputfeedback:Deaecto:2020,outputfeedback:Deaecto:2023}; e.g., \cite{output:Colaneri:2008} provide conditions to stabilise continuous- and discrete-time switched systems through a feedback controller based on output measurements, relying on Lyapunov-Metzler inequalities. 
As an additional challenge, the output may be only available at discrete sampling instants, which requires
sampled-data control methods, explored in the context of switched systems, e.g., by \cite{sampled:Hetel:2013,sampled-data:Albea:2024}.

We consider the stabilisation of a class of nominally linear switched systems, possibly including an affine switching term, affected by uncertain Lipschitz nonlinearities, under dwell-time constraints, sampled-data measurements, and output feedback (Section~\ref{sec:problem}).
In Section~\ref{sec:proposal}, we propose a novel switching law that relies on a state observer, along with Lyapunov-Metzler inequalities \citep{affine:Russo:2024} accounting for dwell-time constraints. Specifically, we discuss the observer design in Section~\ref{Sec:ObsDesign} and the design of the observer gains in Section~\ref{Sec:GainsDesign}.
In Section~\ref{sec:switchstrat} we introduce our proposed observer-based sampled-data switching law, subject to dwell-time constraints, and we derive time-dependent LMI conditions ensuring global asymptotic stability (or, in the presence of switching affine terms, ultimate boundedness) of the closed-loop system. We also obtain an estimate of the average quadratic cost, together with a bound on its maximum deviation from the actual cost function. Moreover, in the presence of switching affine terms with a known common upper bound, we establish explicit ultimate bounds.
In the case of known affine terms, Section~\ref{Sec:ModifCLaw} presents a modified observer-based switching-law design, aimed at improving the estimates on the obtained ultimate bound for the system.
In Section~\ref{Sec:SuffCond}, we discuss the feasibility of the derived LMIs, providing equivalent reduced-order LMI conditions, and prove that the time dependence of the LMIs can be removed by discretisation on a finite grid.
Numerical examples in Section~\ref{sec:sim} illustrate the efficacy of our proposed method, even at stabilising switched systems with unstable modes and with modes having no stable convex combination, concurrently addressing output feedback, sampled-data and dwell-time constraints. The possible presence of affine terms in the dynamics allows us to address real-life engineering applications, such as power converters \citep{switching:Russo:2022,stabilization:Carneiro:2024}, where affine terms usually model the presence of an energy source acting on the system: in particular, we consider the realistic example of a double buck-boost converter, which we use to showcase our approach and compare it with the approach proposed by \cite{sampled:Hetel:2013}, which deals with sampled-data state-feedback stabilisation, but does not include filter design and/or Lyapunov-Metzler inequalities. In the work by \cite{sampled:Hetel:2013}, the dwell time is equal to the inter-sampling time of the state, and hence dwell-time constraints cannot be accommodated into the switching signal. The numerical results show that our methodology performs at least as well under the same conditions and can actually perform better, due to the inherent decoupling between the sampling sequence in the estimation channel and the switching signal updates.
To the best of our knowledge, no existing contribution integrates all the features we consider within a unified observer-based sampled-data framework. Other works in the literature, related to specific aspects of our problem formulation, are discussed in detail in Section~\ref{sec:conclusions}, which also draws conclusions and discusses future work.

\textbf{Notation.} Given matrix $A$, we denote by $A^\top$ its transpose and by $A_{ij}$ its ($i,j$)-th element. When $A$ is symmetric, its entries below the sub-diagonal are denoted by $\star$.  $\R_{\ge 0}$ and 
$\N$  are the sets of non-negative real numbers and natural numbers, respectively. We denote $\mathbb{N}_0=\mathbb{N}\cup\left\{0 \right\}$. For $\ell \in \mathbb{N}$, $[\ell]$ denotes the set $\{1, \dots, \ell\}$, whereas $\mathcal{M}_{\ell}\subseteq \mathbb{R}^{\ell \times \ell}$ denotes the class of irreducible Metzler matrices such that for all matrices $\Pi \in \mathcal{M}_{\ell}$, it holds that $\Pi_{ij}\geq 0$ for all $i\neq j$ and $\sum_{j=1}^{\ell} \Pi_{ij}=0$ for all $i \in [\ell]$. For a symmetric matrix $X\in \mathbb{R}^{n\times n}$, $X \prec 0$ (respectively, $X\succ 0$) indicates that $X$ is positive (respectively, negative) definite.
Finally, $I_n$ denotes the $n$-dimensional identity matrix and $\left\|\cdot \right\|$ the Euclidean norm. The corresponding induced operator norm on $\mathbb{R}^{n\times n}$ is also denoted by $\left\|\cdot \right\|$.
For matrices $\left\{A_i\right\}_{i=1}^k$ of appropriate size, $\operatorname{col}\left\{A_i\right\}_{i=1}^k$ is the block matrix stacking the $A_i$'s in consecutive rows.

\section{Problem Formulation}\label{sec:problem}

Given the state $x\in \R^n$, measured output $y \in \R^m$ and performance output $z \in \R^p$, consider the switched system 
\begin{equation}
	\label{eq:switched_sys} 
\begin{array}{lll}
		&\dot{x}(t) = A_{\sigma(t)}x(t)+\theta b_{\sigma(t)}+f_{\sigma(t)}(x(t)), \quad x(t_0)=x_0,\vspace{0.1cm} \\
		&y(t) = D_{\sigma(t)}x(s_k), \quad t \in [s_k,s_{k+1}),\vspace{0.1cm} \\
		&z(t) = C x(t),
\end{array}
\end{equation}
with \emph{unknown} initial condition $x(t_0)$. We aim to design the piecewise-constant switching law $\sigma\colon \R_{\ge 0} \to [\ell]$, $\ell \in \mathbb{N}$, subject to dwell-time constraints. Thus, denoting with ${\mathcal{T} = \{t_k\}_{k=0}^{\infty}}$ the increasing sequence of switching instants, we introduce a dwell time $T>0$ such that $t_{k+1}-t_k \ge T$ for all $k\in\N$. 
Moreover, the output is sampled at sampling instants $\mathcal{S}= \{s_k\}_{k=0}^\infty$.
Our results allow for the values of the sampling instants $\mathcal{S}$ to be arbitrary, provided that $0<s_{k+1}-s_k \le h$ for all $k\in\N$, where $h>0$ is a maximum inter-sampling time, and $\lim_{k\to\infty} s_k=\infty$. The bound $h$ will be explicitly obtained from suitable matrix inequalities (see e.g., Theorem~\ref{thm:main} in Section~\ref{sec:switchstrat}). Henceforth we will assume, without loss of generality, that $t_0=s_0=0$. Finally, the matrices $A_{\sigma(t)}\in \lbrace A_j\rbrace_{j\in [\ell]}$, $D_{\sigma(t)}\in \lbrace D_j\rbrace_{j\in [\ell]}$ and $C$ are \emph{known}, whereas the vectors $b_{\sigma(t)}\in \left\{b_j \right\}_{j\in [\ell]}$ are unknown; however, an upper bound $\RK{\bar{b}}$ on their norms is available, meaning that 
\begin{equation}\label{eq:UPPBdBi}
\max_{i\in [\ell]}\left\|b_i \right\|\leq \RK{\bar{b}}.
\end{equation}
Similarly, the functions $f_{\sigma(t)}\in \lbrace f_j\rbrace_{j\in [\ell]}$, with $f_j \colon \R^n \to \R^n$, $j\in [\ell]$, are  \emph{unknown}, and represent unknown and unmodelled state-dependent dynamics. The parameter $\theta$ in \eqref{eq:switched_sys} satisfies $\theta \in \left\{0,1 \right\}$, where for $\theta=0$ we will examine the stability of the obtained closed-loop system, whereas for $\theta=1$ we will study its ultimate boundedness.

We make the following standing assumptions.
\begin{assumption}\label{ass:lipschitz}
	Functions $\left\{f_i \right\}_{i\in [\ell]}$ are globally Lipschitz, i.e., there exist \RK{$ \kappa>0$}, such that, for all $x,\,y \in\mathbb{R}^n$,
	\begin{equation}\label{eq:LipKappa}
		\Vert f_i(x)-f_i(y)\Vert\leq \RK{\kappa}\Vert x-y\Vert, \qquad \forall i \in [\ell].
	\end{equation}
	Moreover, for all $i\in [\ell]$, $f_i(0) = 0$. The generality of this assumption is discussed in \RK{Remarks~\ref{rem:LocalLip} and~\ref{Rem:GeneralityAssump}}.
\end{assumption}
\begin{remark}\label{rem_:kappai}
\RK{The subsequent analysis shows that the upper bound $\kappa$ in \eqref{eq:LipKappa} can be replaced by mode-dependent constants $\kappa_i>0$, $i\in [\ell]$. Whenever such tighter mode-dependent bounds are available, the approach presented in this manuscript is expected to yield less conservative results.  }
\end{remark}

\begin{assumption}\label{ass:symultaneous_LF}
	For all $i\in [\ell]$, the pair $(A_i,D_i)$ is observable. Moreover, there exist $\left\{L_i\right\}_{i\in [\ell]}$ such that for all $i\in[\ell]$
	\begin{equation}
        \label{eq:Ui}
		U_i := A_i-L_i D_i
	\end{equation}
	is Hurwitz and the matrices $\lbrace U_i\rbrace_{i\in [\ell]}$ are simultaneously quadratically Lyapunov stable, i.e., there exist $\Omega\in \mathbb{R}^{n \times n}$, with $\Omega\succ0$, and $\eta>0$ such that
	\begin{equation}
        \label{eq:Phi_i}
		\Phi_i(\eta,\Omega) := U_i^\top \Omega+\Omega U_i+2\eta \Omega \prec0, \qquad \forall i \in [\ell].
	\end{equation}
\end{assumption}

We further define the average quadratic cost functional
\begin{equation}
	\label{eq:cost}
	J(t) = \frac{1}{t}\int_{0}^{t} z(\tau)^\top z(\tau) \mathrm{d}\tau,\quad t>0,
\end{equation}
which accounts for the average energy of the performance output $z$.  The value of $J$ may be considered as a measure of the efficiency of the proposed switching algorithm and is of great importance when selecting the best approach in applications, as it serves as a robustness indicator of the process behaviour. Under our assumptions $x(t)$, $t\geq 0$, is not observed, whence $J$ \empty{cannot} be computed directly.

Considering the system \eqref{eq:switched_sys} and the presented standing assumptions, we now formulate our control objectives.
\begin{objective}\label{obj:Objective1}
Using the available sampled-data measurements, design an observer $\varphi(t)$ and an observer-based switching law $\sigma(\varphi(t),t) \colon \R^m \times \R \to [\ell]$ such that the closed-loop system $(x(t),\varphi(t))$ is globally asymptotically stable ($\theta = 0$) or ultimately bounded ($\theta = 1$). In both cases, we require explicit bounds on the decay rate of the closed-loop system and the maximum inter-sampling time $h>0$. For $\theta=1$, we further require a computable upper estimate on the ultimate bound of the state $\left(x(t),\varphi(t) \right)$.
\end{objective}
\begin{objective}\label{obj:Objective2}
Let $\theta = 0$ and assume that an a-priori norm estimate $\left\|x(0)\right\|\leq R$, with $R>0$, on the unknown initial condition is available. Provide a computable estimate $\hat{J}(t)$ of the average cost \eqref{eq:cost} and derive an explicit error bound on $|J(t)-\hat{J}(t)|$, which holds pointwise for all $t\geq 0$. Given $\varepsilon>0$, compute a time $t_*=t_*(\varepsilon,R)>0$ such that $t\geq t_*$ implies $|J(t)-\hat{J}(t) |<\varepsilon$. For $\theta =1$, provide a computable ultimate bound on the discrepancy $\left|J(t)-\hat{J}(t) \right|$ between the actual average cost functional $J(t)$ and the provided estimate $\hat{J}(t)$. 
\end{objective}

\begin{remark}\label{rem:gains}
Assumption~\ref{ass:symultaneous_LF} has been investigated in the literature in the context of polytopic-type uncertainty for parameter-dependent systems of the form $\dot{w} = F(\lambda)w$,
where $\lambda$ is a parameter  (see, e.g.,  \citealp{simulateneousLF:boyd:1989}). Since $\lbrace A_i\rbrace_{i\in [\ell]}$ and $\lbrace D_i\rbrace_{i\in [\ell]}$ are known, Assumption~\ref{ass:symultaneous_LF} can be verified computationally by solving
	\begin{equation*}
		\begin{cases}
			\Omega A_i - Y_i D_i + A_i^\top \Omega - D_i^\top Y_i^{\top} + 2 \eta \Omega\prec 0, \quad i\in [\ell], \\
			\Omega \succ 0
		\end{cases}
	\end{equation*}
	for some $\Omega\succ0$ and $\left\{Y_i\right\}_{i\in [\ell]}$ of appropriate dimensions. If the above inequalities are feasible, $\left\{L_i\right\}_{i\in [\ell]}$ can be set as
	\begin{equation}
		L_i = \Omega^{-1} Y_i, \quad i\in [\ell].
	\end{equation}
By \citep{boyd1994linear} and Finsler's lemma  \citep{Schweppe1973}, Assumption~\ref{ass:symultaneous_LF} holds iff there exist some $\Omega \succ 0$, $\eta>0$ and $\varpi_i\in \mathbb{R},\ i\in [\ell]$, such that for each $i\in [\ell]$ 
\begin{equation*}
\left(A_i+\eta I_n\right)^{\top} \Omega+\Omega\left( A_i+\eta I_n\right)x\prec \varpi_i D_i^{\top}D_i,
\end{equation*}
iff the systems $\dot{z}_i(t)=\left(A_i+\eta I_n \right)z_i(t),\ v_i(t)=D_iz(t)$, $i\in [\ell]$ are jointly dissipative with a common storage function $S(z)=z^{\top}\Omega z$ and supply rates $\omega_i(v_i) = \varpi_iv_i^\top v_i$ \citep{van2000l2}. Section~\ref{Sec:GainsDesign} further discusses the design of the observer gains $\left\{L_i \right\}_{i\in [\ell]}$.  
\end{remark}

\begin{remark}\label{rem:LocalLip}
Our proposed approach is also suitable for the regional stabilisation of systems of the form \eqref{eq:switched_sys}, where $\theta =0$ and the nonlinearities $\left\{f_i \right\}_{i\in [\ell]}$ satisfy $f_i(0)=0$, as well as a \empty{local} Lipschitz condition 
$\left\|f_i(x)-f_i(y) \right\|\leq \RK{\kappa}\left(\rho \right)\left\|x-y \right\|$, for $\|x\|,\|y\|\le \rho$,
with known local Lipschitz constant $\RK{\kappa}(\rho)$ satisfying 
\begin{equation}\label{eq:CondLim}
  \lim_{\rho \to 0^+}\RK{\kappa}(\rho) = 0.  
\end{equation}
Condition \eqref{eq:CondLim} guarantees that, in a neighbourhood of the origin, the linear parts of the dynamics in \eqref{eq:switched_sys} with $\theta=0$ dominate the nonlinearities $\left\{f_i \right\}_{i\in [\ell]}$, thereby making our approach suitable for local stabilisation of \eqref{eq:switched_sys}. Indeed, one can show that subject to \eqref{eq:CondLim},
Theorem~\RK{\ref{thm:main} below} and its proof can be applied locally around the origin. \GG{When mode-dependent Lipschitz constants $\kappa_i(\rho)$ are available (see Remark~\ref{rem_:kappai}), condition \eqref{eq:CondLim} becomes $\lim_{\rho \to 0^+} \max_{i \in [\ell]} \kappa_i(\rho) = 0$.} For concreteness of the presentation, we proceed with the case where $\left\{f_i \right\}_{i\in [\ell]}$ are globally Lipschitz.
\end{remark}

Previously proposed switching strategies such as  those by \cite{practical:sanchez:2019,time-triggered:albea:2021,affine:Russo:2024} cannot be employed to control system \eqref{eq:switched_sys}, since  both $x(0)$ and the right-hand side of \eqref{eq:switched_sys} are unknown.

\section{Observer-based switching law design}\label{sec:proposal}
To achieve the control Objectives~\ref{obj:Objective1} and~\ref{obj:Objective2}, we propose an observer-based sampled-data switching strategy, subject to dwell-time constraints. The observer is designed to compensate the unavailability of $x(t)$, $t\geq 0$ for the design of the switching strategy. For simplicity, we denote the observer-based switching strategy $\sigma(\varphi(t),t)$ as $\sigma(t)$.

\subsection{Observer design}\label{Sec:ObsDesign}
Given any switching law $\sigma(t)$ for system \eqref{eq:switched_sys},
we consider the observer $\varphi \colon \R_{\ge 0} \to \R^n$ satisfying
\begin{equation}\label{eq:obsv}
\begin{array}{lll}
&\dot{\varphi}(t) = A_{\sigma(t)}\varphi(t) + L_{\sigma(t)}[y(s_k)-D_{\sigma(t)}\varphi(s_k)],\vspace{0.1cm}\\
&\hspace{10mm} t\in [s_k,s_{k+1}),\vspace{0.1cm}\\
&\varphi(0)=0,
\end{array}
\end{equation}
where $\left\{L_i \right\}_{i\in [\ell]}$ are observer gains satisfying Assumption~\ref{ass:symultaneous_LF}.
Observer \eqref{eq:obsv} relies solely on knowledge of $\sigma(t)$, $t\geq 0$ (to be designed shortly), of $A_{\sigma(t)}$ and $D_{\sigma(t)}$,  and of the sampled output $y(s_k)=D_{\sigma(t)}x(s_k)$, $t\in [s_k,s_{k+1})$, all of which are assumed to be known. In particular, $b_i$, $i\in [\ell]$, and $f_i$, $i\in [\ell]$, which are assumed unknown, do not appear in the observer design \eqref{eq:obsv}. 

Denote the observer error as $e(t):=x(t)-\varphi(t)$. By \eqref{eq:switched_sys} and \eqref{eq:obsv}, we obtain
\begin{equation*}
\begin{array}{lll}
&\dot{e}(t)=A_{\sigma(t)}e(t)-L_{\sigma(t)}D_{\sigma(t)}e(s_k)+\theta b_{\sigma(t)}+f_{\sigma(t)}(x(t)),\vspace{0.1cm}\\
&\hspace{10mm} t\in [s_k,s_{k+1}] ,\vspace{0.1cm}\\
&e(0)= x(0).
\end{array}
\end{equation*}
Denoting by $\delta_e(t):=e(t)-e(s_k)$, $t\in [s_k,s_{k+1})$, we obtain 
\begin{equation}\label{eq:erroreq1}
\begin{array}{lll}
\dot{e}(t) & = [A_{\sigma(t)}-L_{\sigma(t)}D_{\sigma(t)}]e(t) + L_{\sigma(t)}D_{\sigma(t)}\delta_e(t)\vspace{0.1cm}\\
&\hspace{5mm}+\theta b_{\sigma(t)} + f_{\sigma(t)}(x(t))\vspace{0.1cm} \\
& = U_{\sigma(t)}e(t) + L_{\sigma(t)}D_{\sigma(t)}\delta_e(t)+\theta b_{\sigma(t)} + f_{\sigma(t)}(x(t)),\vspace{0.1cm} \\
&\hspace{5mm}t\in [s_k,s_{k+1}).
\end{array}
\end{equation}
Since $\lim_{k\to \infty}s_k =\lim_{k\to \infty}t_k= \infty$, solutions $\varphi(t)$ and $e(t)$ of \eqref{eq:obsv} and \eqref{eq:erroreq1}, respectively, exist and are absolutely continuous, with the corresponding ODEs satisfied outside of a discrete set of instants $\mathfrak{T}$, such that $\mathfrak{T}$ has a finite intersection with every bounded set in $[0,\infty)$.

Finally, since $x(s_k)-\varphi(s_k)=e(s_k)=e(t)-\delta_e(t)$ and $x(t) = \varphi(t)+e(t)$, we present the closed-loop system in terms of $(\varphi(t),e(t))$ as follows 
\begin{equation}\label{eq:ClosedLoop}
\begin{array}{lll}
\dot{\varphi}(t) \!= A_{\sigma(t)}\varphi(t)+L_{\sigma(t)}D_{\sigma(t)}e(t) - L_{\sigma(t)}D_{\sigma(t)}\delta_e(t) \vspace{0.1cm}\\
\dot{e}(t) \!= U_{\sigma(t)}e(t)+L_{\sigma(t)}D_{\sigma(t)}\delta_e(t)\\
\hspace{30mm}\!+\theta b_{\sigma(t)}+\!f_{\sigma(t)}(\varphi(t)\!+\!e(t)),\ t\geq 0.
\end{array}
\end{equation}
We now proceed to design the switching law $\sigma(t)$, $t\geq 0$, that achieves the control Objectives~\ref{obj:Objective1} and~\ref{obj:Objective2}.

\subsection{Proposed switching strategy}\label{sec:switchstrat}
We first introduce auxiliary variables that we will use in formulating the switching law $\sigma(t)$, $t\geq 0$, and our main result. Consider a scalar $\zeta>0$, matrices $\left\{X_i\right\}_{i\in [\ell]}\subseteq \mathbb{R}^{n\times n}$, $X_i\succ 0$, $i\in [\ell]$, and an irreducible Metzler matrix $\Pi \in \mathcal{M}_{\ell}$ that satisfy the Lyapunov-Metzler inequalities
\begin{equation}\label{eq:LyapMetz}
\begin{array}{lll}
		&A_i^\top X_i + X_i A_i + \sum_{j\in [\ell]\setminus \left\{i \right\}} \Pi_{ij}\left(Y_{1,j}+Y_{2,j}-X_i\right) \\
		&\hspace{35mm}+C^\top C + 2\zeta X_i \prec 0,
\end{array}
	\end{equation}
where, for each $j\in [\ell]$,
\begin{equation}\label{eq:Y12J}
\begin{array}{lll}
& Y_{1,j} = \ee^{(A_j^\top+\zeta I) T}X_j \ee^{(A_j+\zeta I)T},\\
& Y_{2,j} = \int_0^T \ee^{(A_j^\top+\zeta I) \tau}C^\top C\ee^{(A_j+\zeta I)\tau} \mathrm{d}\tau.
\end{array}
\end{equation}
We introduce an observer-based switching law, with dwell-time $T>0$, as follows. Given $t_k\in [0,\infty)$,
	\begin{subequations}
		\label{eq:switching_law}
		\begin{align}
			\sigma(t) & = i, \,\,  t \in [t_k, t_k+T], \label{eq:switching_law_1} \\
			\sigma(t) & = i,  \,\, t >t_k+T, \textrm{ as long as, } \forall j \in [\ell]\setminus \left\{ i\right\}, \label{eq:switching_law_2} \\
			& \textrm{it is } \varphi(t)^\top  (Y_{1,j}+Y_{2,j}) \varphi(t) \geq  \varphi(t)^\top X_i \varphi(t),  \nonumber \\
			\sigma(t_{k+1}) &= \argmin_{j\in [\ell]\setminus \left\{ i\right\}}\left[\varphi(t_{k+1})^\top (Y_{1,j}+Y_{2,j}) \varphi(t_{k+1}) \right],  \label{eq:switching_law_3}
		\end{align}
	\end{subequations}
	where
	\begin{equation*}
		\!\! t_{k+1}:=\!\! \inf_{t > t_k+T}\! \left\lbrace t\ | \ \exists j \colon \varphi(t)^\top \left[ Y_{1,j}\!+\!Y_{2,j}-X_i\right]\varphi(t)\!<\!0 \right\rbrace.
	\end{equation*}
\begin{remark}
As can be seen from \eqref{eq:switching_law}, the set $\left\{t_k \right\}_{k=0}^{\infty}$ can exhibit one of two possible behaviours. First, there may be $t_k$ such that for $t\geq t_k$, the condition in \eqref{eq:switching_law_2} is never violated, and the system ``stabilises'' on a single operational mode: $\sigma(t)\equiv i$ for $t\geq t_k$. In this case, we formally define $t_{k+j} \equiv \infty$ for all $j\geq 1$, thereby guaranteeing that $\lim_{k\to \infty}t_k=\infty$, in the extended real line $\bar{\mathbb{R}}=\mathbb{R}\cup \left\{\infty \right\}$. Second, there may be infinitely many violations of the condition in \eqref{eq:switching_law_2}, thereby leading to $t_k<\infty$ for all $k$ and, since $t_{k+1}-t_k\geq T>0$, $\lim_{k\to \infty}t_k=\infty$ in $\mathbb{R}$. For our stability/ultimate boundedness results, as seen in, e.g., the proof of Theorem~\ref{thm:main}, the distinction between the two behaviours of $\left\{ t_k\right\}_{k=0}^{\infty}$ is immaterial. 
\end{remark}

\begin{remark}
\label{rem:known_bi}
The switching law \eqref{eq:switching_law} requires knowledge of the matrices $A_i$, $i\in [\ell]$, but not of $b_i$, $i\in [\ell]$, and $f_i$, $i\in [\ell]$, which are assumed unknown. When $b_i$, $i\in [\ell]$, are known, incorporating them into both the observer and the switching law designs greatly improves the conservatism of the obtained ultimate bound: in Section~\ref{Sec:ModifCLaw}, we modify our observer-based switching law accordingly.
\end{remark}

We are now ready to state our main result.
\begin{theorem}
    \label{thm:main}
    Consider the closed-loop system \eqref{eq:ClosedLoop} under Assumptions~\ref{ass:lipschitz} and~\ref{ass:symultaneous_LF}, with gains $L_i$ fixed (e.g., as in Remark~\ref{rem:gains}). Let $\zeta,T>0$, $\left\{X_i\right\}_{i\in [\ell]}\subseteq \mathbb{R}^{n \times n}$, $X_i\succ 0$, $i\in [\ell]$, and $\Pi \in \mathcal{M}_{\ell}$ satisfy the Lyapunov-Metzler inequality \eqref{eq:LyapMetz}, and consider the time-varying symmetric positive-definite matrices $P_i\colon \R_{\ge 0} \to \R^{n\times n}$ such that
	\begin{subequations}
		\label{eq:P_dot}
		\begin{align}
			&-\dot{P}_{i}(t) = A_{i}^\top P_{i}(t) + P_{i}(t)A_{i} + C^\top C + 2\zeta P_{i}(t),\nonumber \vspace{0.1cm} \\  & \hspace{10mm}\quad   t \in [t_k, t_k+T),\,  i \in [\ell], \label{eq:Pdot_dyn}\vspace{0.1cm}\\
			&P_{i}(t) = X_{i},\quad  t \in [t_k+T,t_{k+1}), \ i \in [\ell]. \label{eq:Pdot_fin}
		\end{align}
	\end{subequations}
\begin{figure*}[ht]
	\begin{equation}
		\label{eq:Psi}
		\hspace{-1cm}
		\Psi_i(t) := 
		\begin{bsmallmatrix}
			-2(\zeta-\alpha)P_{i}(t) + \gamma \RK{\kappa}^2 I_n & P_{i}(t)L_{i}D_{i}+\gamma \RK{\kappa}^2 I_n & -P_{i}(t) L_{i} D_{i} & 0  \\[6pt]
			\star & \Phi_{i}(\alpha,Q)+\gamma\RK{\kappa}^2 I_n + h^2 \ee^{2\alpha h}U_{i}^\top WU_{i} & Q L_{i} D_{i} + h^2 \ee^{2 \alpha h} U_{i}^\top  W L_{i} D_{i} & Q+h^2 \ee^{2 \alpha h} U_{i}^\top  W \\[6pt]
			\star & \star & -\frac{\pi^2}{4}W + h^2 \ee^{2\alpha h} D_{i}^\top  L_{i}^\top  W L_{i} D_{i} & h^2 \ee^{2\alpha h} D_{i}^\top  L_{i}^\top  W \\[6pt]
			\star & \star & \star & -\gamma I_n + h^2 \ee^{2 \alpha h}W
		\end{bsmallmatrix}
	\end{equation}
\end{figure*} 
    Given tuning parameters $h>0$, $\alpha>0$ and $\RK{\kappa}>0$, $i\in [\ell]$, let there exist matrices $Q\succ0$, $W\succ 0$ and scalar $\gamma>0$ such that the LMIs $\Psi_i(t)\prec 0$, with $\left\{\Psi_i(t) \right\}_{i\in [\ell]}$ given in \eqref{eq:Psi}, hold for all $i\in [\ell]$ and all $t\in [0,T]$.
    Then, 
    \begin{enumerate}
        \item[I.] For $\theta =0$ in \eqref{eq:switched_sys}, the observer-based switching law \eqref{eq:switching_law} yields global asymptotic stability of the origin for 
	    the closed-loop system \eqref{eq:ClosedLoop}, meaning that 
        \begin{equation}\label{eq:ExpStab}
       \begin{array}{lll}
       \left\|\varphi(t) \right\|^2+\left\|e(t) \right\|^2\leq Me^{-2\alpha t}\left\|x(0) \right\|^2,\quad t\geq 0,
        \end{array}
        \end{equation}
        for some $M\geq 1$, thereby satisfying Objective~\ref{obj:Objective1}. The estimate $\hat{J}(t)$ of the average cost, given by 
        \begin{equation}
        \label{eq:cost_estimate}
        \hat{J}(t) = \frac{1}{t}\int_{0}^t \varphi(s)^\top C^\top C\varphi(s) \mathrm{d}s, \quad t>0,
        \end{equation}
        satisfies Objective~\ref{obj:Objective2}.
        \item[II.] For $\theta =1$ in \eqref{eq:switched_sys}, the observer-based switching law \eqref{eq:switching_law}
	    guarantees that the solutions of the closed-loop system \eqref{eq:ClosedLoop} are ultimately bounded, meaning that  
        \begin{equation}\label{eq:UltBdFinal}
    \begin{array}{ll}
    &\limsup_{t \to \infty} \left(\left\|\varphi(t) \right\|^2+\left\|e(t) \right\|^2\right) \leq \tfrac{\Upsilon\RK{\bar{b}}^2}{2 \alpha \min\left(\mathcal{C}_1,\mathcal{C}_3 \right)}=:\overline{B},\\
    &\limsup_{t \to \infty} \| x(t) \|^2 \leq 2\overline{B},
    \end{array}
    \end{equation}
    with explicitly computable $\overline{B}$. Moreover, the estimate $\hat{J}(t)$ of the average cost satisfies
     \begin{equation}\label{eq:AsympDiffJ}
    \limsup_{t\to \infty}\left|J(t)-\hat{J}(t) \right|\leq 3\left\| C^{\top}C\right\|\overline{B}^2.
     \end{equation}
\end{enumerate}
\end{theorem}
\begin{remark}\label{rem:periodicity}
Note that \eqref{eq:P_dot} is an \emph{autonomous} system that is solved on each sub-interval $[t_k,t_k+T]$ with the \emph{same} boundary value $P_i(t)=X_i$, for the fixed matrix $X_i$, whereas on $[t_k+T,t_{k+1})$ we have $P_i\equiv X_i$. Thus, it can be seen that $P_i:\mathbb{R}_{\geq 0}\to \mathbb{R}^{n\times n}$ is fully determined by its values on $[0,T]$ (recall that $t_0=0$ by assumption). Hence, $\Psi_i(t)\prec0$ for $t\in [0,T]$ guarantees $\Psi_i(t)\prec 0$ for all $t\geq 0$. Therefore, in Theorem~\ref{thm:main}, $\Psi_i(t) \prec 0$  is required to hold only for $t \in [0, T]$, and not for all $t \in [t_k, t_k + T]$, $k\in \mathbb{N}_0$.
\end{remark}
\begin{pf}
Let $Q,W\in \mathbb{R}^{n\times n}$ with $Q,W\succ0$ and $\alpha 
\in (0,\zeta)$. Consider the candidate Lyapunov functional
\begin{equation}
\label{eq:V}
V(t) = V_\varphi(t)+V_e(t)
\end{equation}
where 
\begin{equation}
\label{eq:V_phi}
V_\varphi(t) = \varphi(t)^\top P_{\sigma(t)}(t)\varphi(t),\quad t\geq 0 
\end{equation}
and 
\begin{equation}\label{eq:Ve}
\begin{array}{lll}
&\hspace{-3mm}V_e(t) = e(t)^\top Qe(t)+ h^2 \ee^{2\alpha h} \int_{s_k}^{t}\ee^{-2\alpha (t-s)}\dot{e}(s)^\top W \dot{e}(s)  \mathrm{d}s \\
&\hspace{5mm}- \frac{\pi^2}{4} \int_{s_k}^{t} \ee^{-2\alpha(t-s)} \delta_e(s)^\top W \delta_e(s) \mathrm{d}s, \quad t \in [s_k,s_{k+1}).
\end{array}
\end{equation}
$V_{\varphi}(t)$ is continuous on each interval $[t_k,t_{k+1})$, $k\in \mathbb{N}_0$, differentiable on each interval $(t_k,t_{k+1})$, $k\in \mathbb{N}_0$ and satisfies
\begin{equation}\label{eq:Vphi_Equiv}
\mathcal{C}_1\left\|\varphi(t) \right\|^2\leq V_{\varphi}(t)\leq \mathcal{C}_2\left\|\varphi(t) \right\|^2,\quad t\geq 0,
\end{equation}
for some explicitly computable $\mathcal{C}_i>0$, $i=1,2$, in view of the positive-definiteness and continuity of $P_i(t)$, $i\in [\ell]$, satisfying \eqref{eq:P_dot}. Due to the exponential Wirtinger inequality \citep{fridman2014introduction} and the fact that $s \mapsto \delta_e(s)= e(s)-e(s_k)$ is absolutely continuous with square-integrable derivative, the functional $V_e(t)$ satisfies 
\begin{equation}\label{eq:Ve_Equiv}
\begin{array}{lll}
&V_e(t) \geq e(t)^{\top}Qe(t)\geq \mathcal{C}_3\left\|e(t) \right\|^2,\quad t\geq 0,\\
&V_e(s_k) = e(s_k)^{\top}Qe(s_k) \leq \mathcal{C}_4 \left\|e(s_k)\right\|^2,\quad k\in \mathbb{N}_0
\end{array}
\end{equation}
for some explicitly computable $\mathcal{C}_i>0$, $i=3,4$. While 
$V_e(t)$ may exhibit jumps on $\lbrace s_k \rbrace_{k=1}^\infty$,  the exponential Wirtinger inequality further guarantees
\begin{equation}\label{eq:Ve_Jump}
\lim_{t\to s_k^-} V_e(t) \ge V_e(s_k).
\end{equation}
Also, since $t_{k+1}-t_k\ge T>0$, the interval $[s_k,s_{k+1})$ intersects $\lbrace t_k \rbrace_{k=0}^\infty$ at a finite number of points, where $\dot{e}$ exhibits a jump discontinuity, due to switching in $\sigma(t)$. Thus, $\dot{e}$ is piecewise continuous on $[s_k, s_{k+1})$ and the function 
\begin{equation*}
    t\mapsto \int_{s_k}^{t}\ee^{-2\alpha(t-s)} \dot{e}(s)^\top W\dot{e}(s)  \mathrm{d}s
\end{equation*}
is absolutely continuous with the derivative existing outside of a discrete set $\mathfrak{T}_1$ that has a finite intersection with any bounded set in $[0,\infty)$. Therefore, $V(t)$ is indeed a candidate Lyapunov functional, which is continuously differentiable outside the discrete set $\mathfrak{T}\cup\mathfrak{T}_1$, whose intersection with every bounded set in $[0,\infty)$ is finite.

We first consider $t$ in the interval $[t_k,t_k+T)$, where we denote $\sigma(t)=i$. Almost everywhere in $[t_k,t_k+T)$, it holds
\begin{equation*}
\begin{array}{lll}
&\hspace{-2mm}\dot{V}_\varphi (t)  = 2\dot{\varphi}(t)^\top P_{i}(t)\varphi(t) + \varphi(t)^\top \dot{P}_{i}(t)\varphi(t) \vspace{0.1cm}\\
&\hspace{3mm} = 2\left[A_{i}\varphi(t) + L_{i}D_{i}(e(t)-\delta_e(t))\right]^\top P_{i}(t)\varphi(t) \vspace{0.1cm}\\
&\hspace{7mm} + \varphi(t)^\top \dot{P}_{i}(t)\varphi(t)\vspace{0.1cm} \\
&\hspace{3mm} = \varphi(t)^\top \left[A_{i}^\top P_{i}(t)+P_{i}(t)A_{i}+\dot{P}_{i}(t)\right]\varphi(t)  \\
&\hspace{7mm} + 2\varphi(t)^\top P_{i}(t)\left[L_{i}D_{i}e(t)-L_{i}D_{i}\delta_e(t)\right]\vspace{0.1cm} \\
&\hspace{3mm} = -\varphi(t)^\top C^\top C\varphi(t) - 2\zeta\varphi(t)^\top P_{i}(t)\varphi(t)\vspace{0.1cm}  \\
&\hspace{7mm} + 2\varphi(t)^\top P_{i}(t)\left[L_{i}D_{i}e(t)-L_{i}D_{i}\delta_e(t)\right].
\end{array}
\end{equation*}
Therefore, 
\begin{equation}\label{eq:V_dot_phi}
\begin{array}{lll}
\dot{V}_\varphi (t) + &2\alpha V_\varphi(t) + \varphi(t)^\top C^\top C\varphi(t) \vspace{0.1cm} \\
&= -2(\zeta-\alpha)\varphi(t)^\top P_{i}(t)\varphi(t)\vspace{0.1cm} \\
&\quad + 2\varphi(t)^\top P_{i}(t)\left[L_{i}D_{i}e(t)-L_{i} D_{i} \delta_e(t)\right].
\end{array}
\end{equation}
Differentiating $V_e(t)$ a.e. in $[t_k,t_k+T)$ yields
\begin{equation*}
\begin{array}{lll}
&\dot{V}_e(t)+2\alpha V_e(t) \vspace{0.1cm}  \\
&\hspace{12mm}= 2e(t)^\top Q\dot{e}(t) + h^2 \ee^{2\alpha h}\dot{e}(t)^\top W\dot{e}(t) \vspace{0.1cm} \\
&\hspace{17mm} + 2\alpha e(t)^{\top}Qe(t) - \frac{\pi^2}{4} \delta_e(t)^\top W\delta_e(t) \vspace{0.1cm} \\
&\hspace{12mm}= 2 e(t)^\top Q\left[U_{i}e(t)+L_{i}D_{i}\delta_e(t)+f_{i}(\varphi(t)+e(t))\right] \vspace{0.1cm} \\
&\hspace{17mm}+2\theta e(t)^{\top}Qb_i+ 2\alpha e(t)^\top Qe(t) \vspace{0.1cm}\\
&\hspace{17mm}+ h^2 \ee^{2\alpha h} \dot{e}(t)^\top W \dot{e}(t)- \frac{\pi^2}{4}\delta_e(t)^\top W\delta_e(t).
\end{array}
\end{equation*}
As a consequence,
\begin{equation}\label{eq:V_dot_e}
\begin{array}{lll}
&\dot{V}_e(t)+2\alpha V_e(t) \vspace{0.1cm} \\
&\hspace{10mm} = e(t)^\top \Phi_{i}(\alpha,Q)e(t)+2\theta e(t)^{\top}Qb_i\vspace{0.1cm} \\
&\hspace{15mm} + 2 e(t)^\top Q\left[L_{i}D_{i}\delta_e(t)+f_{i}(\varphi(t)+e(t))\right]\vspace{0.1cm} \\
&\hspace{15mm} + h^2 \ee^{2\alpha h} \dot{e}(t)^\top W\dot{e}(t) - \frac{\pi^2}{4}\delta_e(t)^\top W\delta_e(t).	
\end{array}
\end{equation}

For simplicity of notation, let us further write 
\begin{equation}\label{eq:ChiDef}
\chi_i(t):=f_i(x(t))=f_{i}\left(\varphi(t)+e(t) \right).
\end{equation}
Then, in view of Assumption~\ref{ass:lipschitz}, 
\begin{equation}\label{eq:Ni_Term}
\begin{array}{lll}
N_{i}(t) &:= \RK{\kappa}^2 \begin{bmatrix}
\varphi(t) \\
e(t)
\end{bmatrix}^\top 
\begin{bmatrix}
I_n & I_n \\
I_n & I_n
\end{bmatrix}
\begin{bmatrix}
\varphi(t) \\
e(t)
\end{bmatrix}- \Vert \chi_i(t)\Vert^2 \vspace{0.1cm}\\
& = \RK{\kappa}^2 \left\|\varphi(t)+e(t) \right\|^2- \Vert \chi_i(t)\Vert^2 \vspace{0.1cm}\\
& = \RK{\kappa}^2 \left\|x(t) \right\|^2- \Vert \chi_i(t)\Vert^2\geq 0,\quad t\geq 0.
\end{array}
\end{equation}

Finally, employing \eqref{eq:ClosedLoop}, we have
\begin{equation}\label{eq:e_dot}
\begin{array}{lll}
&\hspace{-2mm} \dot{e}(t)^\top W\dot{e}(t) = \left[U_{i}e(t)+L_{i}D_{i}\delta_e(t)+\theta b_i+\chi_{i}(t)\right]^\top W \vspace{0.1cm} \\
&\hspace{23mm} \times \left[U_{i}e(t)+L_{i}D_{i}\delta_e(t)+\theta b_i+\chi_{i}(t)\right] \vspace{0.1cm}  \\
&\hspace{6mm} = \begin{bmatrix}
e(t) \\
\delta_e(t) \\
\chi_{i}(t)\\
b_i
\end{bmatrix}^\top 
\begin{bmatrix}
U_{i}^\top  \\
D_{i}^\top L_{i}^\top \\
I_n\\
\theta I_n
\end{bmatrix}
W
\begin{bmatrix}
U_{i}^\top  \\
D_{i}^\top L_{i}^\top \\
I_n\\
\theta I_n
\end{bmatrix}^\top 
\begin{bmatrix}
e(t) \\
\delta_e(t) \\
\chi_{i}(t)\\
b_i
\end{bmatrix}.
\end{array}
\end{equation}

Recall \eqref{eq:UPPBdBi} and denote
\begin{equation}\label{eq:xidef}
\begin{array}{lll}
&\xi(t) = \mathrm{col}\lbrace \varphi(t),\, e(t),\, \delta_e(t),\, \chi_{i}(t), b_i \rbrace.
\end{array}
\end{equation}
Combining \eqref{eq:V_dot_phi}, \eqref{eq:V_dot_e}, \eqref{eq:Ni_Term} and \eqref{eq:e_dot}, and letting $\gamma,\Upsilon>0$, a.e. in $ [t_k,t_{k}+T)$ we have 
\begin{equation*}
\begin{array}{lll}
&\dot{V}(t)+2\alpha V(t)+ \varphi(t)^\top C^\top C\varphi(t)-\theta \Upsilon \RK{\bar{b}}^2 \vspace{0.1cm} \\
&\hspace{3mm}\leq \dot{V}(t)+2\alpha V(t)+ \varphi(t)^\top C^\top C\varphi(t)-\theta \Upsilon \left\|b_i \right\|^2  +\gamma N_{i}(t) \vspace{0.1cm}\\
&\hspace{3mm} = \xi(t)^\top  \widetilde{\Psi}_{i}(t)\xi(t),
\end{array}
\end{equation*}
where 
\begin{equation}\label{eq:PsiTilde}
\begin{array}{lll}
\small \widetilde{\Psi}_i(t) = \left[\begin{array}{c|c}
\Psi_i(t) & \begin{matrix}0 \\
\theta \left(Q+h^2\ee^{2\alpha h}U_i^{\top}W \right)\\
\theta h^2\ee^{2\alpha h}D_i^{\top}L
_i^{\top}W\\
\theta h^2\ee^{2\alpha h}W\end{matrix}\\ \hline
* & -\theta \left(\Upsilon I_n+h^{2}e^{2\alpha h}W \right)
\end{array}
\right]
\end{array}, \quad t\in[t_k,t_k+T)
\end{equation}
with $\Psi_i(t)$ given in \eqref{eq:Psi}. By Schur complement and taking $\Upsilon\to \infty$, it is clear that $\Psi_i(t)\prec 0$ for all $t\in [0,T]$ guarantees the existence of a large enough $\Upsilon>0$ such that $\widetilde{\Psi}(t)\prec 0$ for all $t\in [t_k,t_k+T)$. Thus, $\Psi_{i}(t)\prec 0$ for all   $t \in [0,T]$ guarantees, for some $\Upsilon>0$, that  
\begin{equation*}
\begin{array}{lll}
&\dot{V}(t)+2\alpha V(t)+ \varphi(t)^\top C^\top C\varphi(t)-\theta \Upsilon \RK{\bar{b}}^2  \le 0 \vspace{0.1cm}\\
&\hspace{40mm}  \text{a.e. in }[t_k,t_k+T).
\end{array}
\end{equation*}

Consider now the interval $[t_k+T, t_{k+1})$, where $\sigma(t)=i$ by \eqref{eq:switching_law}, as the switching condition has not been activated. Taking the derivative of $V(t)$ a.e. in $[t_k+T,t_{k+1})$, we have that $\dot{V}_e(t)$ remains as in \eqref{eq:V_dot_e}. On the other hand, the derivative of $V_\varphi(t) = \varphi(t)^\top X_{i}\varphi(t)$, yields
\begin{equation}\label{eq:SecondInterv}
\begin{array}{lll}
\dot{V}_\varphi(t) & = 2\varphi(t)^\top X_{i}\dot{\varphi}(t) \\
& = 2\varphi(t)^\top X_{i}(A_{i}\varphi(t)+L_{i} D_{i} e(t) - L_{i} D_{i}\delta_e(t))\\
& = \varphi(t)^\top (X_{i} A_{i} + A_{i}^\top  X_{i})\varphi(t) + 2\varphi(t)^\top  X_{i} L_{i} D _{i} e(t) \\
& \quad - 2\varphi(t)^\top X_{i} L_{i} D_{i} \delta_e(t).
\end{array}
\end{equation}
Since no switching has occurred on $[t_k+T,t_{k+1})$, \eqref{eq:switching_law_2} implies that, for all $j\in [\ell]\setminus \left\{i \right\}$,
\begin{equation*}
 \varphi(t)^\top (Y_{1,j}+Y_{2,j}-X_{i})\varphi(t) \ge 0,\quad t\in [t_k+T,t_{k+1}),   
\end{equation*}
whence, taking into account that $\Pi\in \mathcal{M}_{\ell}$,
\begin{equation*}
\begin{array}{lll}
\dot{V}_\varphi(t) &\le  \varphi(t)^\top \bigg(X_{i} A_{i} + A_{i}^\top  X_{i}   \\
&  \quad + \sum_{j\in[\ell]\setminus \left\{i \right\}} \Pi_{ij} (Y_{1,j}+Y_{2,j}-X_{i}) \bigg)\varphi(t) \vspace{0.1cm}\\
&\quad + 2 \varphi(t)^\top X_{i} L_{i} D_{i} e(t) - 2\varphi(t)^\top  X_{i} L_{i} D_{i} \delta_e(t)\vspace{0.1cm} \\
& \leq -\varphi(t)^\top C^\top C\varphi(t) - 2\zeta \varphi(t)^\top X_{i}\varphi(t) \vspace{0.1cm}\\
& \quad +2\varphi(t)^\top X_{i} L_{i} D_{i} e(t) - 2\varphi(t)^\top  X_{i} L_{i} D_{i} \delta_e(t).
\end{array}
\end{equation*}
As a consequence,
\begin{equation}\label{eq:Second_Int_Vphi_Dot}
\begin{array}{lll}
&\dot{V}_\varphi(t)+2\alpha V_\varphi(t) + \varphi(t)^\top C^\top C\varphi(t)\vspace{0.1cm} \\
& \leq -2(\zeta-\alpha) \varphi(t)^\top X_{i}\varphi(t) + 2\varphi(t)^\top X_{i} L_{i} D_{i} e(t)\vspace{0.1cm} \\
& \quad - 2\varphi(t)^\top  X_{i} L_{i} D_{i} \delta_e(t)\quad  \text{a.e. in } [t_k+T,t_{k+1}).
\end{array}
\end{equation}
Recalling \eqref{eq:Pdot_fin}, the upper bound in \eqref{eq:Second_Int_Vphi_Dot} can be obtained by replacing $P_i(t)$ in the right-hand side of \eqref{eq:V_dot_phi} with $P_i(T)=X_i$.
Taking into account \eqref{eq:Second_Int_Vphi_Dot} and applying arguments similar to \eqref{eq:V_dot_e}-\eqref{eq:xidef}, we have
\begin{equation}
\begin{array}{lll}
&\dot{V}(t)+2\alpha V(t) \vspace{0.1cm}\\
&\hspace{8mm}+ \varphi(t)^\top C^\top C\varphi(t) -\theta \Upsilon \RK{\bar{b}}^2 \le \xi(t)^\top \widetilde{\Psi}_i(T)\xi(t)
\end{array}
\end{equation}
a.e. in $[t_k+T,t_{k+1})$, where we employ $P_i(T)=X_i$. Thus, $\Psi_i(t)\prec 0$, $t\in[0,T]$, with $\Psi_i(t)$ given in \eqref{eq:Psi}, guarantees 
\begin{equation*}
\begin{array}{lll}
&\dot{V}(t)+2\alpha V(t)+ \varphi(t)^\top C^\top C\varphi(t)-\theta \Upsilon \RK{\bar{b}}^2 \le 0,\\
&\hspace{40mm}  \text{a.e. in }[t_k+T,t_{k+1}).
\end{array}
\end{equation*}

We now show that $V(t)$ does not exhibit upward jumps at the instants $\left\{t_k \right\}_{k=1}^{\infty}$. By \eqref{eq:V}, since
$V_e(t)$ exhibits only downward jumps, whenever it is discontinuous, as per \eqref{eq:Ve_Jump}, it suffices to show that $V_{\varphi}(t)$ does not exhibit upward jumps on $\left\{t_k \right\}_{k=1}^{\infty}$ (recall that $t_0=0$, by assumption). Let $k\in \mathbb{N}_0$ and assume that $\sigma(t_{k+1}^-)=i$ and $\sigma(t_{k+1})=j$. Due to continuity of $\varphi(t)$, it holds that
\begin{equation*}
\begin{array}{lll}
V_\varphi(t_{k+1}^-) & = \varphi(t_{k+1})^\top X_{i}\varphi(t_{k+1}),\vspace{0.1cm}\\
V_\varphi(t_{k+1}) & = \varphi(t_{k+1})^\top P_{j}(t_{k+1})\varphi(t_{k+1}).
\end{array}
\end{equation*}
By integrating \eqref{eq:P_dot}, we have $P_{j}(t_{k+1})=Y_{1,j}+Y_{2,j}$. Hence, 
\begin{equation*}
V_\varphi(t_{k+1})-V_\varphi(t_{k+1}^-)=\varphi(t_{k+1})^\top [Y_{1,j}+Y_{2,j}-X_{i}]\varphi(t_{k+1}).
\end{equation*}
Recalling the switching condition \eqref{eq:switching_law_3}, at $t=t_{k+1}$ there exists $j_*\in [\ell]\setminus \left\{i \right\}$ such that
\begin{equation*}
\varphi(t_{k+1})^\top [Y_{1,j_*}+Y_{2,j_*}-X_{i}]\varphi(t_{k+1}) \leq 0
\end{equation*} 
and 
\begin{equation*}
\sigma(t_{k+1}) =  \argmin_{j\in [\ell]\setminus\left\{i \right\}} 	\varphi(t_{k+1})^\top [Y_{1,j}+Y_{2,j}]\varphi(t_{k+1}).
\end{equation*}
Therefore, recalling that $\sigma(t_{k+1})=j$, 
\begin{equation*}
\begin{array}{lll}
\varphi(t_{k+1})^\top [Y_{1,j}+Y_{2,j}]\varphi(t_{k+1})\vspace{0.1cm}\\ \le \varphi(t_{k+1})^\top [Y_{1,j_*}+Y_{2,j_*}]\varphi(t_{k+1}) \le \varphi(t_{k+1})^\top X_{i}\varphi(t_{k+1}),
\end{array}
\end{equation*}
implying 
\begin{equation}
\label{eq:delta_V}
V_\varphi(t_{k+1})-V_\varphi(t_{k+1}^-) \leq 0,
\end{equation}
whence also $V_\varphi$ cannot exhibit upward jumps. 

Summarising, we have that 
\begin{equation}
\label{eq:Vdot_ineq}
\begin{array}{lll}
&\dot{V}(t)+2\alpha V(t) + \varphi(t)^\top C^\top C\varphi(t)-\theta \Upsilon \RK{\bar{b}}^2\le 0,\vspace{0.1cm}\\
&\hspace{40mm}\text{ a.e. in } [t_k,t_{k+1}),
\end{array}
\end{equation}
and functional $V(t)$ exhibits no upward jumps at $\left\{t_k \right\}_{k=1}^{\infty}$. Thus, for the case $\theta = 0$,
multiplying \eqref{eq:Vdot_ineq} by $\ee^{2\alpha t}$ and integrating it from $t_k$ to $t$, with $t\in [t_k,t_{k+1})$, yields
\begin{equation*}
\ee^{2\alpha t}V(t)-\ee^{2\alpha t_k}V(t_k)+\int_{t_k}^t\ee^{2\alpha s}\varphi(s)^{\top}C^{\top}C\varphi(s)\mathrm{d}s \leq 0.
\end{equation*}
Summing the latter inequality over $[t_{k'},t_{k'+1}]$, with $k'=0,\dots,k-1$ and $t_0=0$, and $[t_k,t]$ yields
\begin{equation}\label{eq:Vintegration}
\begin{array}{lll}
&\ee^{2\alpha t}V(t)+\int_0^t\ee^{2\alpha s}\varphi(s)^{\top}C^{\top}C\varphi(s)\mathrm{d}s\leq V(0).
\end{array}
\end{equation}
In particular, we obtain $V(t) \le \ee^{-2\alpha t}V(0)$ for all $t\geq 0$. Recalling \eqref{eq:Vphi_Equiv} and \eqref{eq:Ve_Equiv}, and employing $\varphi(0)=0$ and $e(0)=x(0)$, the exponential decay of $V(t)$ yields \eqref{eq:ExpStab} with $M:=\frac{\mathcal{C}_4}{\min\left(\mathcal{C}_1,\mathcal{C}_3 \right)}$, which is explicitly computable. Control Objective~\ref{obj:Objective1} is thus fulfilled for $\theta =0$. Furthermore, consider the average cost estimate $\hat{J}(t)$ defined in \eqref{eq:cost_estimate}, which is bounded on $[0,\infty)$.
Indeed, \eqref{eq:Vintegration} implies that $\int_0^t \varphi(s)^{\top}C^{\top}C\varphi(s)\mathrm{d}s\leq V(0)$; then,
recalling that $V(0)=x(0)^{\top}Qx(0)$ and taking $t\to \infty$, we have that $\sup_{t\geq 0}t\hat{J}(t)\leq x(0)^{\top}Qx(0)$.
Assume now that $\left\|x(0) \right\|\leq R$ for some $R>0$. Then, letting $\tau>0$, we have
\begin{equation}\label{eq:J_equality}
\begin{array}{lll}
J(\tau) & = \frac{1}{\tau} \int_{0}^{\tau} x(s)^\top C^\top Cx(s) \,\mathrm{d}s \vspace{0.1cm} \\
& = \frac{1}{\tau} \int_{0}^{\tau} (\varphi(s)^\top +e(s)^\top )C^\top C(\varphi(s)+e(s)) \,\mathrm{d}s   \\
& = \hat{J}(\tau)  + \frac{2}{\tau} \int_{0}^{\tau} e(s)^\top C^\top C\varphi(s) \,\mathrm{d}s  \\
& \quad + \frac{1}{\tau} \int_{0}^{\tau} e(s)^\top C^\top Ce(s) \,\mathrm{d}s.
\end{array}
\end{equation}
Taking into account \eqref{eq:ExpStab}, we have
\begin{equation*}
\begin{array}{lll}
|J(\tau)-\hat{J}(\tau)| &\leq \tfrac{3M\left\|C^{\top}C \right\|}{2\alpha \tau}\left\|x(0) \right\|^2(1-\ee^{-2\alpha \tau})\vspace{0.1cm}\\
&\leq \tfrac{3M\left\|C^{\top}C \right\|R^2}{2\alpha \tau},
\end{array}
\end{equation*}
where the constants in the final upper bound are computable explicitly. Hence, given $\varepsilon>0$, one can compute a time $t_*=t_*(\varepsilon,R)$ such that $t\geq t_*$ implies $| J(t)-\hat{J}(t)|<\varepsilon$. Hence, $\hat{J}(t)$ fulfils control Objective~\ref{obj:Objective2}. 

Similarly, for the case $\theta = 1$, \eqref{eq:Vdot_ineq} yields
\begin{equation}\label{eq:ineqbound}  
\dot V(t)+2\alpha V(t)+\varphi(t)^\top C^\top C \varphi(t) - \Upsilon \RK{\bar{b}}^2 \leq 0,
\end{equation}
whence arguments similar to those leading to \eqref{eq:Vintegration} yield 
\begin{equation}\label{eq:UltimateBd}
V(t) \leq \ee^{-2 \alpha t} V(0) + \frac{\Upsilon \RK{\bar{b}}^2}{2 \alpha}.   
\end{equation}
Using \eqref{eq:Vphi_Equiv} and \eqref{eq:Ve_Equiv},  we obtain, for $t\geq 0$,
\begin{equation*}
\left\|\varphi(t) \right\|^2+\left\|e(t) \right\|^2\leq Me^{-2\alpha t}\left\|x(0) \right\|^2 + \tfrac{\Upsilon \RK{\bar{b}}^2}{2 \alpha \min\left(\mathcal{C}_1,\mathcal{C}_3 \right)},
\end{equation*}
with $M:=\frac{\mathcal{C}_4}{\min\left(\mathcal{C}_1,\mathcal{C}_3 \right)}$.
Since $\|x(t)\|^2 \leq 2 \| \varphi(t) \|^2 + 2 \| e(t) \|^2$, the explicit ultimate bounds \eqref{eq:UltBdFinal} follow. Control Objective~\ref{obj:Objective1} is thus fulfilled for $\theta=1$. Furthermore, since both $\left\|e(t) \right\|$ and $\left\|\varphi(t) \right\|$ are bounded via $\overline{B}+o(1)$, $t\to \infty$, with $\overline{B}$ given in \eqref{eq:UltBdFinal}, 
employing \eqref{eq:J_equality} finally yields \eqref{eq:AsympDiffJ}.
$\hfill\blacksquare$
\end{pf}
\RK{
\begin{remark}
If one is interested in the stabilisation or ultimate-boundedness of affine switched systems, subject to dwell time, the inequality \eqref{eq:Psi} with $\kappa=0$ provides sufficient conditions to achieve this control objective. Indeed, since $f_i(0)=0$ for all $i\in [\ell]$, in view of Assumption~\ref{ass:lipschitz}, setting the Lipschitz constant $\kappa=0$ is equivalent to setting $f_i\equiv 0$ for all $i\in [\ell]$. Similarly, setting $h=0$ in \eqref{eq:Psi} can be shown to recover the case of continuous-time measurements of the output $y(t)$, $t\geq 0$.
\end{remark}
}
\begin{remark}\label{Rem:GeneralityAssump}
Theorem~\ref{thm:main} can also be applied to the design of an observer-based switching strategy, subject to minimum dwell-time, for the system 
\begin{equation}\label{eq:LipSystem}
	\begin{cases}
		\dot{x}(t) = A_{\sigma(t)}x(t)+g_{\sigma(t)}(x(t)), \quad x(t_0)=x_0, \\
		y(t) = D_{\sigma(t)}x(s_k), \quad t \in [s_k,s_{k+1}),
	\end{cases}
\end{equation}
where $g_{\sigma(t)}\in \left\{g_j \right\}_{j\in [\ell]}$ and $g_j$, $j\in [\ell]$, are globally Lipschitz. Indeed, system \eqref{eq:LipSystem} can be presented as \eqref{eq:switched_sys} with $\theta=1$, $b_{\sigma(t)}=g_{\sigma(t)}(0)$ and $f_{\sigma(t)}(x) = g_{\sigma(t)}(x)-g_{\sigma(t)}(0)$, with the latter being globally Lipschitz and equal to zero at the origin. If $g_j(0)$, $j\in [\ell]$, are unknown, but $\max_{i\in [\ell]}\left\|g_i(0) \right\|$ can be upper bounded, Theorem~\ref{thm:main} can be applied to obtain ultimate boundedness of the closed-loop system \eqref{eq:LipSystem}, with an explicit estimate of the ultimate bound. If $g_j(0)$, $j\in [\ell]$, are known, we refer the reader to the modified switching law design of Section~\ref{Sec:ModifCLaw}.
\end{remark}

\subsection{Design of the observer gains $\left\{L_i \right\}_{i\in [\ell]}$}\label{Sec:GainsDesign}
Theorem~\ref{thm:main} assumes that the gains $L_i$, $i\in [\ell]$, are given and proceeds to obtain $Q\succ 0$ such that $\Phi_i(\alpha,Q)\prec 0$ $\forall i\in [\ell]$; see the $(2,2)$ block in \eqref{eq:Psi}. The design of the observer gains  $\left\{L_i \right\}_{i\in [\ell]}$ via the approach in Remark~\ref{rem:gains} may lead to significant conservatism in verifying $\Psi_i(t)\prec 0,\ i\in [\ell], \ t\in [0,T]$, with $\Psi_i(t)$ given in \eqref{eq:Psi}. Indeed, such a a design of $\left\{L_i \right\}_{i\in [\ell]}$ does not take into account that the obtained gains are to be substituted into $\Psi_i(t)\prec 0,\ i\in [\ell], \ t\in [0,T]$ as exogenous parameters and thereby does not factor the presence of the other decision variables $Q$ and $W$, as well as the matrices $X_i$, coming from \eqref{eq:LyapMetz}, in the inequalities to be verified.

Thus, it is of interest to design the observer gains from the matrix inequalities obtained in Theorem~\ref{thm:main}. The following proposition, where we assume for simplicity of presentation that $f_i\equiv 0$ for all $i\in [\ell]$, provides \emph{sufficient} conditions for the design of $L_i,\ i\in [\ell]$ via Bilinear Matrix Inequalities (BMIs). BMIs\AR{, known to be $\mathcal{NP}$-hard \citep{BMI:Toker:1995},} can be solved using, e.g., Alternating Convex Search and Augmented Lagrangian Methods (see \citealp{boyd2004convex} and \citealp{vanantwerp2000tutorial}), which are implemented in existing software packages, including PENBMI \AR{(see \citealp{PENBMI:Kocvara:2003} and references therein)} and BMISOLVER \AR{(see \citealp{BMIsolver:Dinh:2012,BMIsolver:Dinh:2012b} and references therein)}.

\begin{prop}\label{Prop:LiDesign}
Let $\zeta,T>0$, $\left\{X_i\right\}_{i\in [\ell]}\subseteq \mathbb{R}^{n \times n}$, $X_i\succ 0$, $i\in [\ell]$, and $\Pi \in \mathcal{M}_{\ell}$ satisfy the Lyapunov-Metzler inequality \eqref{eq:LyapMetz}, and consider the time-varying symmetric positive-definite matrices $P_i\colon \R_{\ge 0} \to \R^{n\times n}$ satisfying \eqref{eq:P_dot}. Given tuning parameters $h>0$, $\alpha>0$, let there exist matrices $Q\succ 0$, $W\succ 0$, $R\succ 0$, $Y_i\in \mathbb{R}^{n\times p}$, $i\in [\ell]$, and $J\in \mathbb{R}^{n\times 2n}$ such that 
\begin{equation}\label{eq:BMIs}
\begin{array}{lll}
 &\begin{bmatrix}
    R & 0\\
    0& -W
\end{bmatrix}+\begin{bmatrix}
    I_n\\
    -W
\end{bmatrix}J+J^{\top}\begin{bmatrix}
    I_n & -W
\end{bmatrix} \prec 0,\vspace{0.2cm}\\
&\Lambda_i(t) \prec 0,\quad i\in [\ell], \ t\in [0,T],
\end{array}
\end{equation}
where $\left\{\Lambda_i(t) \right\}_{i\in [\ell]}$ are given in \eqref{eq:Psi1R}. Then, setting the observer gains as $L_i = Q^{-1}Y_i$, the inequalities $\Psi_i(t)\prec 0,\ i\in [\ell], \ t\in [0,T]$, with $\left\{\Psi_i(t) \right\}_{i\in [\ell]}$ given in \eqref{eq:Psi}, hold with $L_i,\ i\in [\ell]$ and the same values of $\zeta$, $\alpha$, $h$, $W$, $Q$.
\end{prop}
\begin{remark}
The inequalities $\Lambda_i(t)\prec 0$, $i\in [\ell]$, $t\in [0,T]$ are linear in $Q$, $Y_i$, $i\in [\ell]$, $W$ and $R$. Thus, the only source of non-linearity in \eqref{eq:BMIs} arises from the first inequality, which is quadratic in $W$ and $J$. Furthermore, as the proof of Proposition~\ref{Prop:LiDesign} shows, only the matrix $W\succ 0$ prevents the conversion of $\Psi_i(t)\prec 0,\ i\in [\ell], \ t\in[0,T]$ into LMIs for the design of the gains. Recall that $W$ emanates from the component of $V_e(t)$ in \eqref{eq:Ve}, which is aimed at compensating the sampling in the state measurements.    
\end{remark}

\begin{pf}
Existence of $J\in \mathbb{R}^{n\times 2n}$ such that the first inequality in \eqref{eq:BMIs} holds 
implies $W^{-1}\succ R$, i.e., $\vartheta^{\top}W^{-1}\vartheta> \vartheta^{\top}R\vartheta$ for any $\vartheta\in \mathbb{R}^{n}\setminus \left\{0 \right\}$. In fact, $W^{-1}\succ R$ is equivalent to 
\begin{equation}\label{eq:finsler1}
\begin{array}{lll}
\zeta^{\top}H^{\top}H\zeta = 0, \ \zeta \neq 0 \Longrightarrow \zeta^{\top}G\zeta<0,
\end{array}
\end{equation}
where
\begin{equation*}
\begin{array}{lll}
&H = \begin{bmatrix}
I & -W
\end{bmatrix}, \ G = \begin{bmatrix}
R & 0 \\
0 & -W
\end{bmatrix},\ \zeta =\begin{bmatrix}
\vartheta\\ z   
\end{bmatrix}.
\end{array}
\end{equation*}
By Finsler's lemma \citep{Schweppe1973}, \eqref{eq:finsler1} is equivalent to the existence of $\mu \in \mathbb{R}$ such that $G-\mu H^{\top}H\prec 0$. Employing \citep[Section 2.6.2]{boyd1994linear}, the latter holds iff there exists $J\in \mathbb{R}^{n\times 2n}$ such that 
\begin{equation*}
\begin{array}{lll}
&G+H^{\top}J+J^{\top}H \prec 0,
\end{array}
\end{equation*}
which is exactly the first inequality in \eqref{eq:BMIs}.

Recall \eqref{eq:e_dot} in the proof of Theorem~\ref{thm:main}, and denote $Y_i = QL_i$. In view of the assumption that $f_i\equiv 0,\ i\in [\ell]$, we can consider $\Psi_i(t)$ in \eqref{eq:Psi} with the last row and the last column deleted, as well as set $\RK{\kappa}= 0$. Performing Schur complement, the inequality $\Psi_i(t)\prec 0,\ t\in [0,T]$  is then equivalent to $\Psi_i^{(1)}\left(t;W^{-1}\right)\prec 0,\ t\in [0,T]$, with $\Psi_i^{(1)}\left(t;W^{-1}\right)$ given in \eqref{eq:Psi1R}. In view of $W^{-1}\succ R$, which has been previously demonstrated, we have $-h^{-2}e^{-2\alpha h}W^{-1}\prec -h^{-2}e^{-2\alpha h}R$, whence $\Psi_i^{(1)}\left(t;W^{-1} \right)\preceq \Psi_i^{(1)}\left(t;R \right),\ t\in [0,T]$, where $\Psi_i^{(1)}\left(t;R \right)$ is obtained from $\Psi_i^{(1)}\left(t;W^{-1} \right)$ by replacing $W^{-1}$ with $R$ in the $(4,4)$ entry. 

Next, denoting 
\begin{equation}\label{eq:ZiDes}
\begin{array}{lll}
Z_1 = \small \begin{bmatrix}
P_i(t)\\0\\0\\0
\end{bmatrix},\ Z_2 = \small \begin{bmatrix}
0 \\ D_i^{\top}Y_i^{\top}\\
-D_i^{\top}Y_i^{\top}\\0
\end{bmatrix},\ Z_3 = -Z_2,\  Z_4 =\small \begin{bmatrix}
  0\\0\\0\\I  
\end{bmatrix},
\end{array}
\end{equation}
we see that 
\begin{equation*}
\begin{array}{lll}
&\Psi_i^{(1)}\left(t;R \right) = \Psi_i^{(2)}\left(t;R \right)\\
&\hspace{20mm}+\sum_{i=1}^2\left(Z_{2i-1}Q^{-1}Z_{2i}^{\top}+Z_{2i}Q^{-1}Z_{2i-1}^{\top}\right),
\end{array}
\end{equation*}
with $\Psi_i^{(2)}\left(t;R \right)$ given in \eqref{eq:Psi1R}. Performing a square completion, for $i=1,2$ we have 
\begin{equation*}
\begin{array}{lll}
&Z_{2i-1}Q^{-1}Z_{2i}^{\top}+Z_{2i}Q^{-1}Z_{2i-1}^{\top}\vspace{0.1cm} \\
&\hspace{5mm}= \frac{1}{2}\left(Z_{2i-1}+Z_{2i} \right)Q^{-1}\left(Z_{2i-1}+Z_{2i} \right)^{\top}\\
&\hspace{10mm}-\frac{1}{2}\left(Z_{2i-1}-Z_{2i} \right)Q^{-1}\left(Z_{2i-1}-Z_{2i} \right)^{\top},
\end{array}
\end{equation*}
with the rightmost term being negative semidefinite. Thus, we obtain
\begin{equation*}
\begin{array}{lll}
&\Psi_i^{(1)}\left(t;R \right) \leq  \Psi_i^{(2)}\left(t;R \right)\\
&\hspace{20mm}+\frac{1}{2}\sum_{i=1}^2\left(Z_{2i-1}+Z_{2i} \right)Q^{-1}\left(Z_{2i-1}+Z_{2i} \right)^{\top}.
\end{array}
\end{equation*}
By applying Schur complement to the right-hand side, we conclude that $\Psi_i^{(1)}\left(t;R \right)\prec 0,\ i\in [\ell], \ t\in [0,T]$ holds if $\Lambda_i(t)\prec 0,\ i\in [\ell], \ t\in [0,T]$, with $\Lambda_i(t)$ given in \eqref{eq:Psi1R}. $\hfill\blacksquare$
\end{pf}
\begin{figure*}[ht]
\begin{equation}
\label{eq:Psi1R}
\begin{array}{lll}
&\Psi_i^{(1)}\left(t;W^{-1}\right) := 
\begin{bsmallmatrix}
-2(\zeta-\alpha)P_{i}(t) & P_{i}(t)Q^{-1}Y_{i}D_{i}& -P_{i}(t)Q^{-1} Y_{i} D_{i} & 0  \\[6pt]
\star & A_i^{\top}Q+QA_i-D_i^{\top}Y_i^{\top}-Y_iD_i+2\alpha Q & Y_{i} D_{i} & A_i^{\top}-D_i^{\top}Y_i^{\top}Q^{-1} \\[6pt]
\star & \star & -\frac{\pi^2}{4}W &  D_i^{\top}Y_i^{\top}Q^{-1} \\[6pt]
\star & \star & \star &  -h^{-2}\ee^{-2\alpha h}W^{-1}
\end{bsmallmatrix},\vspace{0.2cm}\\
&\Psi_i^{(2)}\left(t;R\right) := 
\begin{bsmallmatrix}
-2(\zeta-\alpha)P_{i}(t) & 0& 0 & 0  \\[6pt]
\star & A_i^{\top}Q+QA_i-D_i^{\top}Y_i^{\top}-Y_iD_i+2\alpha Q & Y_{i} D_{i} & A_i^{\top} \\[6pt]
\star & \star & -\frac{\pi^2}{4}W &  0\\[6pt]
\star & \star & \star &  -h^{-2}\ee^{-2\alpha h}R
\end{bsmallmatrix},\ \Lambda_i(t) = \small\left[
\begin{array}{c|c}
\Psi_i^{(2)}\left(t;R \right) & \begin{matrix}
    P_i(t) & 0\\
    D_i^{\top}Y_i^{\top} & -D_i^{\top}Y_i^{\top}\\
    -D_i^{\top}Y_i^{\top} & D_i^{\top}Y_i^{\top}\\
    0 & I
\end{matrix}\\ \hline\star & -2\operatorname{diag}\left(Q,Q \right)
\end{array}.
 \right]
\end{array}
\end{equation}
\end{figure*}

The BMI-based 
design procedure given in Proposition~\ref{Prop:LiDesign} may be computationally challenging. However, the proof of Proposition~\ref{Prop:LiDesign} suggests two ways in which such a design may be relaxed in order to recover LMIs for the design of the observer gains $\left\{L_i \right\}_{i\in [\ell]}$. 

\begin{prop}\label{Prop:QeqW}
Let $\zeta,T>0$, $\left\{X_i\right\}_{i\in [\ell]}\subseteq \mathbb{R}^{n \times n}$, $X_i\succ 0$, $i\in [\ell]$, and $\Pi \in \mathcal{M}_{\ell}$ satisfy the Lyapunov-Metzler inequality \eqref{eq:LyapMetz}, and consider the time-varying symmetric positive-definite matrices $P_i\colon \R_{\ge 0} \to \R^{n\times n}$ satisfying \eqref{eq:P_dot}. Given tuning parameters $h>0$, $\alpha>0$, let there exist matrices $Q\succ 0$ and $Y_i\in \mathbb{R}^{n\times p}$ such that $\Theta_i(t)\prec 0,\ i\in [\ell], t\in[0,T]$, with $\left\{\Theta_i(t) \right\}_{i\in [\ell]}$ given in \eqref{eq:Thetai}. Then, setting the observer gains as $L_i = Q^{-1}Y_i$, the inequalities $\Psi_i(t)\prec 0,\ i\in [\ell], \ t\in [0,T]$, with $\left\{\Psi_i(t) \right\}_{i\in [\ell]}$ given in \eqref{eq:Psi}, hold with $L_i,\ i\in [\ell]$ and the same values of $\zeta$, $\alpha$, $h$, $Q$ and $W=Q$.
\end{prop}
\begin{figure*}[ht]
\begin{equation}
\label{eq:Thetai}
\begin{array}{lll}
&\tilde{\Psi}_i(t) := 
		\begin{bsmallmatrix}
			-2(\zeta-\alpha)P_{i}(t) & P_{i}(t)L_{i}D_{i}& -P_{i}(t) L_{i} D_{i}  \\[6pt]
			\star & \Phi_{i}(\alpha,Q)+ h^2 \ee^{2\alpha h}U_{i}^\top QU_{i} & Q L_{i} D_{i} + h^2 \ee^{2 \alpha h} U_{i}^\top  Q L_{i} D_{i} \\[6pt]
			\star & \star & -\frac{\pi^2}{4}Q + h^2 \ee^{2\alpha h} D_{i}^\top  L_{i}^\top  Q L_{i} D_{i} 
\end{bsmallmatrix}, \ \Theta_i^1(t) = \small \begin{bmatrix}
    P_i(t) & 0  \\
    D_i^{\top}Y_i^{\top} & A_i^{\top}Q-D_i^{\top}Y_i^{\top} \\
    -D_i^{\top}Y_i^{\top} & D_i^{\top}Y_i^{\top} 
\end{bmatrix},\vspace{0.1cm}\\ [6pt]
&\Theta_i^{0}(t) := 
\begin{bsmallmatrix}
-2(\zeta-\alpha)P_{i}(t) & 0& 0   \\[6pt]
\star & A_i^{\top}Q+QA_i-D_i^{\top}Y_i^{\top}-Y_iD_i+2\alpha Q & Y_{i} D_{i} &  \\[6pt]
\star & \star & -\frac{\pi^2}{4}Q 
\end{bsmallmatrix},\ \Theta_i(t) = \small\left[
\begin{array}{c|c}
\Theta_i^{0}(t) & \Theta^1_i(t) \\ \hline\star & \begin{matrix}
-2Q & 0 \\
\star & -2h^{-2}\ee^{-2\alpha h}Q 
\end{matrix}
\end{array}
 \right].
\end{array}
\end{equation}
\end{figure*}
\begin{pf}
The proof follows arguments similar to those of Proposition~\ref{Prop:LiDesign}. Thus, we only provide a sketch.

In view of the assumption that $f_i\equiv 0,\ i\in [\ell]$, we can consider $\Psi_i(t)$ in \eqref{eq:Psi} with the last row and the last column deleted, as well as set $\RK{\kappa}= 0$. Then, setting $W=Q$ and $Y=QL_i$, we note that $\Psi_i(t)$ becomes $\tilde{\Psi}_i(t)$, given in \eqref{eq:Thetai}. We further have
\begin{equation*}
\begin{array}{lll}
&U_i^{\top}QU_i = \left(QA_i-Y_iD_i \right)^{\top}Q^{-1}\left(QA_i-Y_iD_i \right),\vspace{0.1cm}\\
&U_i^{\top}QL_iD_i = (QA_i-Y_iD_i)^{\top}Q^{-1}Y_iD_i,\vspace{0.1cm}\\
&D_i^{\top}L_i^{\top}QL_iD_i = D_{i}^{\top}Y_i^{\top}Q^{-1}Y_iD_i.
\end{array}
\end{equation*}
Thus,
\begin{equation*}
\begin{array}{lll}
&\tilde{\Psi}_i(t)=\Theta_i^{0}(t)+Z_1Q^{-1}Z_2^{\top}+Z_2Q^{-1}Z_1^{\top} \vspace{0.1cm} \\
&\hspace{22mm}+h^{2}e^{2\alpha h}V_1Q^{-1}V_1^{\top},
\end{array}
\end{equation*}
where $Z_1$ and $Z_2$ are as in \eqref{eq:ZiDes}, but with the last entry deleted, and
\begin{equation*}
\begin{array}{lll}
&V_1 = \small \begin{bmatrix}
  0 \\ \left(QA_i-Y_iD_i \right)^{\top}\\
  D_i^{\top}Y_i^{\top}
\end{bmatrix}.
\end{array}
\end{equation*}
The result of the proposition now follows by applying a square completion, followed by Schur complement, similarly to the proof of Proposition~\ref{Prop:LiDesign}.$\hfill\blacksquare$ 
\end{pf}

An alternative design approach, which leads to even simpler LMIs, is as follows. Considering $\Psi_i(t)$ in \eqref{eq:Psi} and set $h=0$, $\alpha=0$, $\RK{\kappa}=0$ and $Y_i=QL_i, \ i\in [\ell]$, it can be easily seen that $\Psi_i(t)\prec 0, i\in [\ell]$ hold iff 
\begin{equation}\label{eq:LMIDesign}
\begin{array}{lll}
& \small \begin{bmatrix}
   -2\zeta P_i(t) & P_i(t)Q^{-1}Y_iD_i\vspace{0.1cm} \\
   \star & A_i^{\top}Q+QA_i-D_i^{\top}Y_i^{\top}-Y_iD_i
\end{bmatrix}\prec 0,\vspace{0.1cm}\\
&\hspace{40mm} i\in[\ell],\ t\in[0,T].
\end{array}
\end{equation}
By performing square completion as in the proof of Proposition~\ref{Prop:LiDesign}, feasibility of \eqref{eq:LMIDesign} is guaranteed if 
\begin{equation}\label{eq:LMIDesign1}
\begin{array}{lll}
&\small \begin{bmatrix}
   -2\zeta P_i(t) & 0 & P_i(t)\\
   \star & A_i^{\top}Q+QA_i-D_i^{\top}Y_i^{\top}-Y_iD_i & D_i^{\top}Y_i^{\top}\\
   \star & \star & -2Q
\end{bmatrix}\prec 0,\vspace{0.1cm}\\
&\hspace{40mm}i\in [\ell], \ t\in [0,T].
\end{array}
\end{equation}
The LMIs \eqref{eq:LMIDesign1} can be used to obtain the gains $L_i = Q^{-1}Y_i,\ i\in [\ell]$. Then, the obtained gains can be employed in verifying $\Psi_i(t)\prec 0,\ i\in [\ell],\ t\in [0,T]$ with $h>0$, $\alpha>0$ and $\RK{\kappa}>0$.

\subsection{Modified switching law for the case of known $b_i,\ i\in [\ell]$\label{Sec:ModifCLaw}}

The analysis of Section~\ref{Sec:ObsDesign} assumes that the affine terms $\left\{b_i \right\}_{i\in [\ell]}$ are unknown and only the bound \eqref{eq:UPPBdBi} is available. As a result, the observer design in \eqref{eq:obsv} and the proposed switching law \eqref{eq:switching_law} do not take into account the presence of $\left\{b_i \right\}_{i\in [\ell]}$ in the dynamics of \eqref{eq:switched_sys} when $\theta=1$. If $\left\{b_i \right\}_{i\in [\ell]}$ are known, simply employing Theorem~\ref{thm:main} in order to obtain an estimate on the ultimate bound may lead to conservative results, due to inefficient tracking of $x(t)$ by $\varphi(t)$ and unnecessary switches of the switching law  \eqref{eq:switching_law}. In this section, we propose a modified switching law tailored to the case where $\left\{b_i \right\}_{i\in [\ell]}$ are known. 

Consider \eqref{eq:switched_sys} with $\theta =1$ and $\left\{b_i \right\}_{i\in [\ell]}$ known and let 
\begin{equation*}
\begin{array}{lll}
&\bar{x}(t) = \left[\begin{array}{c} x(t) \\ \hline 1\end{array}\right]\in \mathbb{R}^{n+1},\ \bar{C} = \left[
\begin{array}{c|c}
C & 0
\end{array}
\right]\in \mathbb{R}^{p\times (n+1)},\vspace{0.2cm} \\
&\bar{A}_{\sigma(t)} = \left[\begin{array}{c|c}
A_{\sigma(t)} & b_{\sigma(t)}\\ \hline
0 & 0
\end{array}\right]\in \mathbb{R}^{(n+1)\times (n+1)}, \vspace{0.2cm} \\
&\bar{f}_{\sigma(t)}\left(\bar{x}(t) \right) = \left[\begin{array}{c}f_{\sigma(t)}\left(x(t) \right)\\ \hline 0  \end{array}\right]\in \mathbb{R}^{n+1} ,\vspace{0.2cm}\\
&\bar{D}_{\sigma(t)} = \left[\begin{array}{c|c}
D_{\sigma(t)} & 0
\end{array}\right]\in \mathbb{R}^{m\times(n+1)}.
\end{array}
\end{equation*}
Then, \eqref{eq:switched_sys} can be embedded into the extended system 
\begin{equation}\label{eq:switched_sys1}
\begin{array}{lll}
		&\dot{\bar{x}}(t) = \bar{A}_{\sigma(t)}\bar{x}(t)+\bar{f}_{\sigma(t)}(\bar{x}(t)), \  \bar{x}(t_0)= \left[ \begin{array}{c} x_0 \\ \hline 1 \end{array}\right],\vspace{0.1cm} \\
		&y(t) = \bar{D}_{\sigma(t)}\bar{x}(s_k), \quad t \in [s_k,s_{k+1}),\vspace{0.1cm} \\
		&z(t) = \bar{C} \bar{x}(t),
\end{array}   
\end{equation}
where $\left\{\bar{A}_i \right\}_{i\in [\ell]}$ are known, due to our assumption on availability of $\left\{b_i \right\}_{i\in [\ell]}$. In view of \eqref{eq:switched_sys1}, let
\begin{equation*}
\begin{array}{lll}
\bar{\varphi}(t) = \left[\begin{array}{c} \varphi(t) \\ \hline 1 \end{array} \right]\in \mathbb{R}^{n+1}, 
\ \bar{L}_{\sigma(t)} = \left[\begin{array}{c} L_{\sigma(t)}\\ \hline 0
\end{array}\right]\in \mathbb{R}^{(n+1)\times m},
\end{array}
\end{equation*}
where the observer $\bar{\varphi}(t)$ satisfies
\begin{equation}\label{eq:obsv1}
\begin{array}{lll}
&\dot{\bar{\varphi}}(t) = \bar{A}_{\sigma(t)}\bar{\varphi}(t)+ \bar{L}_{\sigma(t)}[y(s_k)-\bar{D}_{\sigma(t)}\bar{\varphi}(s_k)],\vspace{0.1cm}\\
&\hspace{30mm} t\in [s_k,s_{k+1}),\vspace{0.1cm}\\
&\bar{\varphi}(0)=\left[\begin{array}{c} 0 \\ \hline 1\end{array} \right].
\end{array}
\end{equation}
Here, $\left\{b_i \right\}_{i\in [\ell]}$ are included in the observer design through $\left\{\bar{A}_i \right\}_{i\in [\ell]}$, whereas the observer gains $\left\{L_i \right\}_{i\in [\ell]}$ satisfy Assumption~\ref{ass:symultaneous_LF}. We denote 
\begin{equation}\label{eq:EstErrMod}
\begin{array}{lll}
&\bar{e}(t)= \bar{x}(t)-\bar{\varphi}(t) = \left[\begin{array}{c} e(t) \\ 0\end{array}\right] = \left[\begin{array}{c} x(t)-\varphi(t) \\ 0 \end{array} \right],\vspace{0.1cm}\\
&\bar{\delta}_e(t) =\left[\begin{array}{c}  \delta_e(t) \\ 0\end{array}\right] =\bar{e}(t)-\bar{e}(s_k) ,\vspace{0.1cm}\\
&\bar{U}_{\sigma(t)} = \bar{A}_{\sigma(t)}-\bar{L}_{\sigma(t)}\bar{D}_{\sigma(t)}.
\end{array}
\end{equation}
Then, by arguments similar to  \eqref{eq:erroreq1}--\eqref{eq:ClosedLoop}, we obtain the following closed-loop system for $\left(\bar{\varphi}(t),\bar{e}(t) \right)$
\begin{equation}\label{eq:ClosedLoopMod}
\begin{array}{lll}
\dot{\bar{\varphi}}(t) = \bar{A}_{\sigma(t)}\bar{\varphi}(t)+\bar{L}_{\sigma(t)}\bar{D}_{\sigma(t)}\bar{e}(t)-\bar{L}_{\sigma(t)}\bar{D}_{\sigma(t)}\bar{\delta}_e(t),\vspace{0.1cm}\\
\dot{\bar{e}}(t) =\bar{U}_{\sigma(t)}\bar{e}(t)+\bar{L}_{\sigma(t)}\bar{D}_{\sigma(t)}\bar{\delta}_e(t)+\bar{f}_{\sigma(t)}\left( \bar{\varphi}(t)+\bar{e}(t)\right).
\end{array}
\end{equation}

For the switching strategy design, consider given $\varepsilon,\zeta>0$, matrices $\left\{X_i \right\}_{i\in [\ell]}\subseteq \mathbb{R}^{n \times n}$, $X_i\succ 0,\ i\in [\ell]$ and an irreducible Metzler matrix $\Pi\in \mathcal{M}_{\ell}$. We further denote
\begin{equation}\label{eq:XEMod}
\begin{array}{lll}
& \bar{X}_i = \left[\begin{array}{c|c} X_i & 0 \\ \hline 0 & 1\end{array} \right],\quad E = \left[\begin{array}{c|c}0 & 0 \\ \hline 0 & 1 \end{array} \right].
\end{array}
\end{equation}
Let the following Lyapunov-Metzler inequalities be satisfied
\begin{equation}\label{eq:LyapMetzMod}
\begin{array}{lll}
		&\bar{A}_i^\top \bar{X}_i + \bar{X}_i \bar{A}_i + \sum_{j\in [\ell]\setminus \left\{i \right\}} \Pi_{ij}\left(\bar{Y}_{1,j}+\bar{Y}_{2,j}-\bar{X}_i\right) \\
		&\hspace{35mm}+\bar{C}^\top \bar{C} + 2\zeta \bar{X}_i \prec \varepsilon E,
\end{array}
	\end{equation}
where $\bar{Y}_{1,j}$ and $\bar{Y}_{2,j}$ are given by the expressions in \eqref{eq:Y12J}, with $A_j$, $X_j$ and $C$ replaced by $\bar{A}_j$, $\bar{X}_j$ and $\bar{C}$, respectively. Then, the observer-based switching law, with dwell time $T>0$, is defined as follows. Given $t_k\in [0,\infty)$,
	\begin{subequations}
		\label{eq:switching_lawMod}
		\begin{align}
			\sigma(t) & = i, \,\,  t \in [t_k, t_k+T], \label{eq:switching_law_1Mod} \\
			\sigma(t) & = i,  \,\, t >t_k+T, \textrm{ as long as, } \forall j \in [\ell]\setminus \left\{ i\right\}, \label{eq:switching_law_2Mod} \\
			& \textrm{it is } \bar{\varphi}(t)^\top  (\bar{Y}_{1,j}+\bar{Y}_{2,j}) \bar{\varphi}(t) \geq  \bar{\varphi}(t)^\top \bar{X}_i \bar{\varphi}(t),  \nonumber \\
			\sigma(t_{k+1}) &= \argmin_{j\in [\ell]\setminus \left\{ i\right\}}\left[\bar{\varphi}(t_{k+1})^\top (\bar{Y}_{1,j}+\bar{Y}_{2,j}) \bar{\varphi}(t_{k+1}) \right],  \label{eq:switching_law_3Mod}
		\end{align}
	\end{subequations}
	where
	\begin{equation*}
		\!\! t_{k+1}:=\!\! \inf_{t > t_k+T}\! \left\lbrace t\ | \ \exists j \colon \bar{\varphi}(t)^\top \left[ \bar{Y}_{1,j}\!+\!\bar{Y}_{2,j}-\bar{X}_i\right]\bar{\varphi}(t)\!<\!0 \right\rbrace.
	\end{equation*}
We now state the main theorem of this section, which is analogous to Theorem~\ref{thm:main} for the case of known $\left\{ b_i\right\}_{i\in [\ell]}$.
\begin{theorem} \label{thm:mainMod}
Consider the closed-loop system \eqref{eq:ClosedLoopMod} under Assumptions~\ref{ass:lipschitz} and~\ref{ass:symultaneous_LF}, with gains $L_i$ fixed (e.g., as in Remark~\ref{rem:gains}). Let $\varepsilon, \zeta,T>0$, $\left\{X_i\right\}_{i\in [\ell]}\subseteq \mathbb{R}^{n \times n}$, $X_i\succ 0$, $i\in [\ell]$, and $\Pi \in \mathcal{M}_{\ell}$ satisfy the Lyapunov-Metzler inequality \eqref{eq:XEMod}, \eqref{eq:LyapMetzMod} and consider the time-varying symmetric positive-definite matrices $\bar{P}_i\colon \R_{\ge 0} \to \R^{n\times n}$ such that
	\begin{subequations}
		\label{eq:P_dotMod}
		\begin{align}
			&-\dot{\bar{P}}_{i}(t) = \bar{A}_{i}^\top \bar{P}_{i}(t) + \bar{P}_{i}(t)\bar{A}_{i} + \bar{C}^\top \bar{C} + 2\zeta \bar{P}_{i}(t),\nonumber \\  & \hspace{10mm}\quad   t \in [t_k, t_k+T),\,  i \in [\ell], \label{eq:Pdot_dynMod}\\
			&\bar{P}_{i}(t) = \bar{X}_{i},\quad  t \in [t_k+T,t_{k+1}), \ i \in [\ell]. \label{eq:Pdot_finMod}
		\end{align}
	\end{subequations}
\begin{figure*}[ht]
	\begin{equation}
		\label{eq:PsiMod}
	\begin{array}{lll}
		&\bar{\Psi}_i(t) := 
		\begin{bsmallmatrix}
			-2(\zeta-\alpha)\bar{P}_{i}(t) + \gamma \RK{\kappa}^2 \Omega_1 & \bar{P}_{i}(t)\bar{L}_{i}D_i+\gamma \RK{\kappa}^2 \Omega_2 & -\bar{P}_{i}(t) \bar{L}_{i} D_i & 0  \\[6pt]
			\star & \Phi_{i}(\alpha,Q)+\gamma\RK{\kappa}^2 I_n + h^2 \ee^{2\alpha h}U_{i}^\top WU_{i} & Q L_{i} D_{i} + h^2 \ee^{2 \alpha h} U_{i}^\top  W L_{i} D_{i} & Q+h^2 \ee^{2 \alpha h} U_{i}^\top  W \\[6pt]
			\star & \star & -\frac{\pi^2}{4}W + h^2 \ee^{2\alpha h} D_{i}^\top  L_{i}^\top  W L_{i} D_{i} & h^2 \ee^{2\alpha h} D_{i}^\top  L_{i}^\top  W \\[6pt]
			\star & \star & \star & -\gamma I_n + h^2 \ee^{2 \alpha h}W
		\end{bsmallmatrix},\vspace{0.1cm} \\
        & \Omega_1 = \small \begin{bmatrix}
    I_n & 0 \\ 0 & 0
\end{bmatrix}\in \mathbb{R}^{(n+1)\times (n+1)},\quad \Omega_2 = \begin{bmatrix}
    I_n\\ 0
\end{bmatrix}\in \mathbb{R}^{(n+1)\times n}
\end{array}
	\end{equation}
\end{figure*} 
Given tuning parameters $h>0$, $\alpha>0$ and $\RK{\kappa}>0$,  let there exist matrices $Q\succ0$, $W\succ 0$ and scalar $\gamma>0$ such that the LMIs $\bar{\Psi}_i(t)\prec 0$ hold for all $i\in [\ell]$ and all $t\in [0,T]$, where $\left\{\bar{\Psi}_i(t) \right\}_{i\in [\ell}$ are given in \eqref{eq:PsiMod}. 
    Then, the observer-based switching law \eqref{eq:switching_lawMod}
	guarantees that the solutions of the closed-loop system \eqref{eq:ClosedLoopMod} satisfy the ultimate bounds
    \begin{equation}\label{eq:UltBdFinalMod}
\begin{array}{ll}
&\limsup_{t \to \infty} \left(\left\|\varphi(t) \right\|^2+\left\|e(t) \right\|^2\right) \leq \tfrac{\varepsilon}{2 \alpha \min\left(\mathcal{C}_1,\mathcal{C}_3 \right)}=:\overline{B}_1,\\
&\limsup_{t \to \infty} \| x(t) \|^2 \leq 2\overline{B}_1,
\end{array}
\end{equation}
with explicitly computable positive constants $\mathcal{C}_1$ and $\mathcal{C}_3$. Furthermore, the estimate $\hat{J}(t)$ of the average cost $J(t)$ satisfies \eqref{eq:AsympDiffJ} with $\overline{B}$ replaced by $\overline{B}_1$.
\end{theorem}
\begin{pf}
The proof follows the same arguments as the proof of Theorem~\ref{thm:main}. Thus, we only elaborate on the main differences. 

Let $Q,W\in \mathbb{R}^{n\times n}$ with $Q,W\succ 0$ and $\alpha\in (0,\zeta)$ and consider the candidate Lyapunov functional 
\begin{equation}
\label{eq:Vmod}
V(t) = V_{\bar{\varphi}}(t)+V_e(t),
\end{equation}
where, in view of the structure of $\bar{e}(t)$ in \eqref{eq:EstErrMod}, $V_e(t)$ is given by \eqref{eq:Ve} and $V_{\bar{\varphi}}(t)$ is given by
\begin{equation}
\label{eq:V_phiMod}
V_{\bar{\varphi}}(t) = \bar{\varphi}(t)^\top \bar{P}_{\sigma(t)}(t)\bar{\varphi}(t),\quad t\geq 0.
\end{equation}

Consider the interval $[t_k,t_k+T)$, where we denote $\sigma(t)=i$. By arguments similar to \eqref{eq:V_dot_phi}--\eqref{eq:V_dot_e}, almost everywhere on $[t_k,t_k+T)$, it holds
\begin{equation}\label{eq:VphiMod}
\begin{array}{lll}
&\dot{V}_{\bar{\varphi}}(t)+2\alpha V_{\bar{\varphi}}(t)+\bar{\varphi}(t)^{\top}\bar{C}^{\top}\bar{C}\bar{\varphi}(t)\vspace{0.1cm}\\
&\hspace{8mm}=\dot{V}_{\bar{\varphi}}(t)+2\alpha V_{\bar{\varphi}}(t)+\varphi(t)^{\top}C^{\top}C\varphi(t)\vspace{0.1cm}\\
&\hspace{8mm}= -2(\zeta-\alpha)\bar{\varphi}(t)^{\top}\bar{P}_i(t)\bar{\varphi}(t)+2\bar{\varphi}(t)^{\top}\bar{P}_i(t)\bar{L}_i\bar{D}_i\bar{e}(t)\vspace{0.1cm}\\
&\hspace{12mm}-2\bar{\varphi}(t)^{\top}\bar{P}_i(t)\bar{L}_i\bar{D}_i\bar{\delta}_e(t)
\end{array}
\end{equation}
and 
\begin{equation}\label{eq:V_dot_eMod}
\begin{array}{lll}
&\dot{V}_e(t)+2\alpha V_e(t)= e(t)^\top \Phi_{i}(\alpha,Q)e(t) \vspace{0.1cm} \\
&\hspace{10mm} +2 e(t)^\top QL_iD_i\delta_e(t)+ 2 e(t)^\top Q\chi_i(t) \vspace{0.1cm} \\
&\hspace{13mm} + h^2 \ee^{2\alpha h} \dot{e}(t)^\top W\dot{e}(t) - \frac{\pi^2}{4}\delta_e(t)^\top W\delta_e(t),	
\end{array}
\end{equation}
where we used the notation \eqref{eq:ChiDef}. Recalling \eqref{eq:PsiMod}, the upper bound \eqref{eq:Ni_Term} is replaced by
\begin{equation}\label{eq:Ni_Term_Mod}
\begin{array}{lll}
&\bar{N}_{i}(t)= \RK{\kappa}^2 \begin{bmatrix}
\bar{\varphi}(t) \\
e(t)
\end{bmatrix}^\top \left[\begin{array}{c|c}
\Omega_1 & \Omega_2\\
\hline \star & I_n
\end{array}\right]
\begin{bmatrix}
\bar{\varphi}(t) \\
e(t)
\end{bmatrix}- \Vert \chi_i(t)\Vert^2 \geq 0,
\end{array}
\end{equation}
whereas \eqref{eq:e_dot} is replaced by 
\begin{equation}\label{eq:e_dotMod}
\begin{array}{lll}
&\hspace{-2mm} \dot{e}(t)^\top W\dot{e}(t) = \begin{bmatrix}
e(t) \\
\delta_e(t) \\
\chi_{i}(t)
\end{bmatrix}^\top 
\begin{bmatrix}
U_{i}^\top  \\
D_{i}^\top L_{i}^\top \\
I_n
\end{bmatrix}
W
\begin{bmatrix}
U_{i}^\top  \\
D_{i}^\top L_{i}^\top \\
I_n
\end{bmatrix}^\top 
\begin{bmatrix}
e(t) \\
\delta_e(t) \\
\chi_{i}(t)
\end{bmatrix}.
\end{array}
\end{equation}
Denoting 
\begin{equation}\label{eq:xidefMod}
\begin{array}{lll}
&\bar{\xi}(t) = \mathrm{col}\lbrace \bar{\varphi}(t),\, e(t),\, \delta_e(t),\, \chi_{i}(t) \rbrace
\end{array}
\end{equation}
and letting $\gamma>0$, a.e. in $[t_k,t_k+T)$, we have 
\begin{equation*}
\begin{array}{lll}
&\dot{V}(t)+2\alpha V(t)+ \varphi(t)^\top C^\top C\varphi(t)\vspace{0.1cm} \\
&\hspace{3mm}\leq \dot{V}(t)+2\alpha V(t)+ \varphi(t)^\top C^\top C\varphi(t) +\gamma \bar{N}_{i}(t) \vspace{0.1cm}\\
&\hspace{3mm} = \bar{\xi}(t)^\top  \bar{\Psi}_{i}(t)\bar{\xi}(t)\leq 0.
\end{array}
\end{equation*}

Considering now the interval $[t_k+T,t_{k+1})$, where $\sigma(t)=i$ by \eqref{eq:switching_lawMod} and the switching condition has not occurred, by arguments similar to \eqref{eq:SecondInterv}--\eqref{eq:Second_Int_Vphi_Dot}, we have  
\begin{equation}\label{eq:VphiMod1}
\begin{array}{lll}
&\dot{V}_{\bar{\varphi}}(t)+2\alpha V_{\bar{\varphi}}(t)+\varphi(t)^{\top}C^{\top}C\varphi(t)\vspace{0.1cm}\\
&\hspace{10mm}= -2(\zeta-\alpha)\bar{\varphi}(t)^{\top}\bar{X}_i\bar{\varphi}(t)+2\bar{\varphi}(t)^{\top}\bar{X}_i\bar{L}_i\bar{D}_i\bar{e}(t)\vspace{0.1cm}\\
&\hspace{14mm}-2\bar{\varphi}(t)^{\top}\bar{X}_i\bar{L}_i\bar{D}_i\bar{\delta}_e(t)+\varepsilon \bar{\varphi}(t)^{\top}E\bar{\varphi}(t)\\
&\hspace{10mm}= -2(\zeta-\alpha)\bar{\varphi}(t)^{\top}\bar{X}_i\bar{\varphi}(t)+2\bar{\varphi}(t)^{\top}\bar{X}_i\bar{L}_i\bar{D}_i\bar{e}(t)\vspace{0.1cm}\\
&\hspace{14mm}-2\bar{\varphi}(t)^{\top}\bar{X}_i\bar{L}_i\bar{D}_i\bar{\delta}_e(t)+\varepsilon,
\end{array}
\end{equation}
where we used the fact that $\bar{\varphi}(t)^{\top}E\bar{\varphi}(t)=1$, in view of \eqref{eq:obsv1} and \eqref{eq:XEMod}.
Letting $\gamma>0$, a.e. in $[t_k+T,t_{k+1})$, we have 
\begin{equation*}
\begin{array}{lll}
&\dot{V}(t)+2\alpha V(t)+ \varphi(t)^\top C^\top C\varphi(t)\vspace{0.1cm} \\
&\hspace{3mm}\leq \dot{V}(t)+2\alpha V(t)+ \varphi(t)^\top C^\top C\varphi(t) +\gamma \bar{N}_{i}(t) \vspace{0.1cm}\\
&\hspace{3mm} = \bar{\xi}(t)^\top  \bar{\Psi}_{i}(T)\bar{\xi}(t)+\varepsilon\leq \varepsilon.
\end{array}
\end{equation*}

Summarising, we have that 
\begin{equation}
\label{eq:Vdot_ineqMod}
\begin{array}{lll}
&\dot{V}(t)+2\alpha V(t) + \varphi(t)^\top C^\top C\varphi(t)\le \varepsilon,\vspace{0.1cm}\\
&\hspace{40mm}\text{ a.e. in } [t_k,t_{k+1}),
\end{array}
\end{equation}
and the functional $V(t)$ exhibits no upward jumps at $\left\{t_k \right\}_{k=1}^{\infty}$. Employing arguments similar to the proof of Theorem~\ref{thm:main}, \eqref{eq:Vdot_ineqMod} yields
\eqref{eq:UltBdFinalMod}, where $\mathcal{C}_1$ and $\mathcal{C}_3$ can be computed as in \eqref{eq:Vphi_Equiv} and \eqref{eq:Ve_Equiv}. The estimate \eqref{eq:AsympDiffJ} with $\overline{B}$ replaced by $\overline{B}_1$, which is given in \eqref{eq:UltBdFinalMod}, follows exactly as in Theorem~\ref{thm:main}. $\hfill\blacksquare$
\end{pf}

\section{LMI feasibility guarantees and discretisation}\label{Sec:SuffCond}
In this section, we study the feasibility of the matrix inequalities obtained in Section~\ref{sec:switchstrat}. For simplicity of presentation, we will focus on Theorem~\ref{thm:main} and the feasibility of $\Psi_i(t)\prec 0$, $i\in [\ell]$, $t\in [0,T]$, defined in \eqref{eq:Psi}. However, we stress that similar arguments can be applied to the design conditions presented in Section~\ref{Sec:GainsDesign} and Section~\ref{Sec:ModifCLaw}.

\begin{proposition}\label{eq:FeasSuff}
Under Assumptions~\ref{ass:lipschitz} and~\ref{ass:symultaneous_LF}, consider $P_i(t)$, $i\in [\ell]$, as in Theorem~\ref{thm:main}. Let $Q\in \mathbb{R}^{n\times n}$ with $Q\succ0$ be such that, for all $i\in [\ell]$,
\begin{equation}\label{eq:GammaDef}
	\Gamma_i(t):=\begin{bsmallmatrix}
		-2 \zeta P_{i}(t)  & P_{i}(t) L_{i} D_{i} \\
		\star & \Phi_{i}(0,Q)
	\end{bsmallmatrix} \prec 0,\quad t\in [0,T].
\end{equation}
Then, $\Psi_i(t)\prec 0$, $t\in[0,T]$, holds for all $i\in [\ell]$, provided that $h,\alpha,\RK{\kappa}>0$, are small enough and $\gamma>0$ is large enough.
\end{proposition}
\begin{pf}
Consider \eqref{eq:Psi} and set $\alpha =0$, $h=\RK{\kappa} = \gamma^{-\beta}$, with $\beta>\frac{1}{2}$, and $W=\gamma I_n$.  Then, the bottom-right $2\times2$ block in \eqref{eq:Psi} has the form
\begin{equation*}
-\gamma \left(\begin{bsmallmatrix}
\frac{\pi^2}{4}I_n& 0\\
0 & I_n
\end{bsmallmatrix} - \gamma^{-2\beta}\begin{bsmallmatrix}
D_i^{\top}L_i^{\top}L_iD_i & D_i^{\top}L_i^{\top}\\
\star & I_n
\end{bsmallmatrix}\right)=:-\gamma \Xi(\gamma).
\end{equation*}
Since $\lim_{\gamma\to \infty}\Xi(\gamma)=\begin{bsmallmatrix} \frac{\pi^2}{4}I_n& 0\\ 0 & I_n \end{bsmallmatrix}\succ 0$,
we have $-\gamma\Xi(\gamma)\prec 0$ for $\gamma\geq \gamma_*>0$,  with $\gamma_*$ large enough, whence $h,\RK{\kappa}>0$ small enough. Then, applying Schur complement in \eqref{eq:Psi}, with respect to this $2\times 2$ block, we have that  $\Psi_i(t)\prec 0$, $t\in [0,T]$, holds if and only if
\begin{equation}\label{eq:SuffCond}
\begin{array}{lll}
\Gamma_i(t) &+ \begin{bsmallmatrix}
\gamma^{1-2\beta}I_n & \gamma^{1-2\beta}I_n\\
\star & \gamma^{1-2\beta}\left(I_n+U_i^{\top}U_i \right)
\end{bsmallmatrix}\vspace{0.1cm}\\
&\hspace{-3mm}+\gamma^{-1}\Delta_i(\gamma,t)\Xi(\gamma)^{-1}\Delta_i(\gamma,t)^{\top}\prec 0, \quad t\in [0,T],
\end{array}
\end{equation}
where we note that $1-2\beta<0$ and define
\begin{equation*}
\Delta_i(\gamma,t): = \begin{bsmallmatrix}
    -P_i(t)L_iD_i & 0\\
    Q L_iD_i+\gamma^{1-2\beta}U_i^{\top}L_i D_i & Q+\gamma^{1-2\beta}U_i^{\top}
\end{bsmallmatrix}.
\end{equation*}
Since $P_i$ is continuous on $[0,T]$, so are $\Gamma_i(t)$ and $\Lambda_i(\gamma,t)$. Thus, in view of compactness of $[0,T]$, there exist some $\mu,\eta>0$ such that, for all $t\in[0,T]$, we have $\Gamma_i(t)\prec -\mu I_n$ and $\sup_{t\in[0,T],\gamma\geq \gamma_*}\left\|\Delta_i(\gamma,t)\right\|\leq \eta$, for all $i\in [\ell]$. Thus, by increasing $\gamma_*$ if necessary, \eqref{eq:SuffCond} can be guaranteed to hold. Taking again into account that $P_i(t)$, $i\in [\ell]$, are continuous on $[0,T]$, similar arguments show there exists some $\alpha_*>0$ such that feasibility of $\Psi_i(t)\prec 0$, $i\in [\ell] $ for all $t\in [0,T]$ is guaranteed for $\alpha\in(0,\alpha_*)$. $\hfill\blacksquare$
\end{pf}
\begin{proposition}\label{prop:FeasNec}
Assume that $\Psi_i(t)\prec 0$, $t\in [0,T]$, for all $i\in [\ell]$, with $\Psi_i(t)$ given in \eqref{eq:Psi}. Then, \eqref{eq:GammaDef} holds for all $i\in [\ell]$ with the same $Q\in\mathbb{R}^{n\times n}$ and $\zeta>0$ as in \eqref{eq:Psi}.
\end{proposition}
\begin{pf}
For $i\in [\ell]$, $\Psi_i(t)\prec 0$, $t\in [0,T]$, implies that 
\begin{equation}\label{eq:22BlockNegative}
\begin{array}{lll}
0\succ& \begin{bsmallmatrix}
-2(\zeta-\alpha)P_i(t)+\gamma \RK{\kappa}^2 I_n & P_i(t)L_iD_i+\gamma \RK{\kappa}^2 I_n\\
\star & \Phi_i(\alpha,Q)+\gamma \RK{\kappa}^2 I_n + h^2e^{2\alpha h}U_i^{\top}W U_i
\end{bsmallmatrix}\vspace{0.1cm}\\
&=\Gamma_i(t) +\begin{bsmallmatrix}
2\alpha P_i(t) & 0\\
0 & 2\alpha Q+h^{2}e^{2\alpha h}U_i^{\top}WU_i
\end{bsmallmatrix}+\gamma \RK{\kappa}^2 \begin{bsmallmatrix}
I_n & I_n \\
I_n & I_n
\end{bsmallmatrix}\vspace{0.1cm}\\
&=:\Gamma_i(t) + M(t) + M_0, \quad t\in [0,T].
\end{array}
\end{equation}
Since $\alpha,h>0$ and $Q,W,P_i(t)\succ0$, $t\in [0,T]$, we have $M(t)\succ 0$, $t\in [0,T]$. Also, since $\gamma,\RK{\kappa}>0$, it follows that $M_0\succ0$; see \eqref{eq:Ni_Term}. Therefore
\begin{equation*}
\begin{array}{lll}
\Gamma_i(t)\prec \Gamma_i(t)+M(t)+M_0 \prec 0,
\end{array}
\end{equation*}
and this concludes the proof. $\hfill\blacksquare$
\end{pf}

Propositions~\ref{eq:FeasSuff} and~\ref{prop:FeasNec} together provide a sufficient and necessary condition for the feasibility of $\Psi_i(t)\prec 0$ for \emph{some} tuning parameters, summarised in the following theorem.

\begin{theorem}
$\Psi_i(t)\prec 0$, $t\in[0,T]$, holds $\forall i\in [\ell]$ for \emph{some} $\alpha>0$, inter-sampling time $h>0$ and Lipschitz constants $\RK{\kappa}>0$, if and only if \eqref{eq:GammaDef} holds $\forall i\in [\ell]$.
\end{theorem}

To verify the feasibility of $\Psi_i(t)\prec 0$, $t\in [0,T]$, with $\Psi_i(t)$ given in \eqref{eq:Psi}, we propose to discretise the interval $[0,T]$ into a densely-spaced grid $0=\tau_0<\dots<\tau_M=T$ and simultaneously verify the LMIs $\Psi_i(\tau_j)\prec 0$, $i\in [\ell]$, $j\in [M]$. The following results formally justify such a discretisation.

\begin{lemma}\label{lem:PiBound}
Let $X_i\in \mathbb{R}^{n\times n}$ with $X_i\succ 0$ be fixed and consider $P_i(t)$ satisfying \eqref{eq:P_dot}. Then,
\begin{equation*}
\max_{t\in [0,T]}\left\|P_i(t) \right\|\leq \left(\left\|X_i \right\|+T\left\|C^{\top}C \right\| \right)\ee^{2T\left\|A_{i,\zeta} \right\|}=: \mathbf{P}_i,
\end{equation*}
where $A_{i,\zeta}:=A_i+\zeta I_n$.
\end{lemma}
\begin{pf}
For $\Sigma_i(t):=P_i(T-t)$, $t \in [0,T]$, in view of \eqref{eq:P_dot}, $\dot \Sigma_i = A_{i,\zeta}^{\top} \Sigma_i + \Sigma_i A_{i,\zeta}+C^\top C$. Integrating on $[0,\tau]$ with $0 \leq \tau \leq T$, since $\Sigma_i(0)=X_i$, yields
\begin{equation}\label{eq:PiIntegral}
\begin{array}{lll}
\Sigma_i(\tau) = X_i+\tau C^{\top}C\\
 \hspace{2mm} +A_{i,\zeta}^{\top}\int_{0}^{\tau}\Sigma_i(s)\mathrm{d}s+\int_{0}^{\tau} \Sigma_i(s)\mathrm{d}s A_{i,\zeta}, \quad \tau \in [0,T].
\end{array}
\end{equation}

Since $\max_{t\in [0,T]}\|P_i(t)\|=\max_{\tau\in [0,T]}\|\Sigma_i(\tau)\|$, the result follows from taking the norm on both sides of \eqref{eq:PiIntegral}, using norm sub-multiplicativity, exchanging the norm with the integrals, and employing the Grönwall inequality.$\hfill\blacksquare$
\end{pf}
\begin{lemma}\label{lem:PiErrorBound}
Let $X_i\in \mathbb{R}^{n\times n}$ with $X_i\succ 0$ be fixed and consider $P_i(t)$ satisfying \eqref{eq:P_dot}. Given $\varepsilon>0$, let $\mu>0$ be such that 
\begin{equation*}
\begin{array}{lll}
\vartheta_i(\mu):=&\mu \ee^{2 \left\|A_{i,\zeta} \right\|\mu}\left\|C^{\top}C \right\|\\
&+\left[2\left(\ee^{\left\|A_{i,\zeta} \right\|\mu}-1 \right)+\left(e^{\left\|A_{i,\zeta} \right\|\mu}-1 \right)^2 \right]\mathbf{P}_i<\varepsilon.
\end{array}
\end{equation*}
Then, for any $0\leq \tau_0\leq \tau_1\leq T$ with $\tau_1-\tau_0<\mu$, it holds $\left\|P_i(\tau_1)-P_i(\tau_0) \right\|<\varepsilon$.
\end{lemma}
\begin{pf}
Employing \eqref{eq:Pdot_dyn}, we have 
\begin{equation*}
\begin{array}{lll}
0 = \frac{\mathrm{d}}{\mathrm{d}t}\left(\ee^{A_{i,\zeta}^{\top}t}P_i(t)\ee^{A_{i,\zeta}t} \right)+\ee^{A_{i,\zeta}^{\top}t}C^{\top}C\ee^{A_{i,\zeta}t}.
\end{array}
\end{equation*}
Integrating the latter equality on $[\tau_0,\tau_1]$, employing some algebraic manipulations and denoting $\Delta \tau:=\tau_1-\tau_0$ and $\mathcal{I}:=\int_0^{\Delta \tau}\ee^{A_{i,\zeta}^{\top}s}C^{\top}C\ee^{A_{i,\zeta}s}\mathrm{d}s$, we end up with
\begin{equation*}
\begin{array}{lll}
\hspace{-2mm}P_i(\tau_0)-P_i(\tau_1)&= \mathcal{I} + \ee^{A_{i,\zeta}^{\top}\Delta\tau}P(\tau_1)\ee^{A_{i,\zeta}\Delta\tau}-P(\tau_1)\\ &= \mathcal{I}+P_i(\tau_1)h_{i}^{\zeta}(\Delta \tau)+h_i^{\zeta}(\Delta \tau)^{\top}P_i(\tau_1)\\& \quad+h_i^{\zeta}(\Delta \tau)^{\top}P_i(\tau_1)h_i^{\zeta}(\Delta \tau),
\end{array}
\end{equation*}
where $h_i^{\zeta}(\mu):=\ee^{A_{i,\zeta}\mu}-I_n$ satisfies $\|h_i^{\zeta}(\mu)\|\leq \ee^{\|A_{i,\zeta}\| \mu}-1$. Taking the norm of $P_i(\tau_1)-P_i(\tau_0)$ and using Lemma~\ref{lem:PiBound}, we obtain
\begin{equation*}
\|P_i(\tau_1)-P_i(\tau_0) \|\leq \vartheta_i(\Delta \tau)\leq \vartheta_i(\mu)<\varepsilon. \hspace{3mm} \blacksquare
\end{equation*}
\end{pf}
\begin{corollary}\label{Cor:PhiDifference}
Consider \eqref{eq:Psi} with fixed scalars $\zeta,h,\gamma,\RK{\kappa}>0$ and $\alpha\in[0,\zeta)$, and fixed matrices $W,Q\in \mathbb{R}^{n\times n}$, $W,Q\succ 0$.
Under the assumptions of Lemma~\ref{lem:PiErrorBound}, let $\mu>0$ be such that $\vartheta_i(\mu)<\frac{\varepsilon}{2\left(\zeta-\alpha+\left\|L_iD_i \right\| \right)}$. Then,
\begin{equation*}
\begin{array}{lll}
\left\|\Psi_i(\tau_1)-\Psi_i(\tau_0) \right\|<\varepsilon.
\end{array}
\end{equation*}
\end{corollary}
\begin{pf}
Employing \eqref{eq:Psi} and Lemma~\ref{lem:PiErrorBound}, we have
\begin{align*}
&\hspace{-2mm}\left\| \Psi_i(\tau_1)-\Psi_i(\tau_0)\right\|\leq 2(\zeta-\alpha+\left\| L_iD_i\right\| )\left\|P_i(\tau_1)-P_i(\tau_0) \right\|\vspace{0.1cm}\\
&\hspace{26mm}\leq 2(\zeta-\alpha+\left\| L_iD_i\right\| )\vartheta_i(\Delta \tau)\\
&\hspace{26mm}\leq 2(\zeta-\alpha+\left\| L_iD_i\right\| )\vartheta_i(\mu) <\varepsilon. \hspace{3mm} \blacksquare
\end{align*}
\end{pf}
We can now formulate the following theorem.
\begin{theorem}\label{Thm:SamplingLMIs}
Let $i\in [\ell]$ be fixed and consider the grid $0=\tau_0<\dots<\tau_M=T$ on $[0,T]
$. Given $\varepsilon,\zeta>0$ and $\alpha\in (0,\zeta)$, assume that 
\begin{equation}\label{eq:Stencil}
 \mathfrak{s}:=\max_{0\leq j \leq M-1}\left(\tau_{j+1}-\tau_j \right)<\mu,   
\end{equation}
with $\mu>0$ such that $\vartheta_i(\mu)<\frac{\varepsilon}{2\left(\zeta-\alpha+\left\|L_i D_i \right\| \right)}$.
Consider $h,\gamma,\RK{\kappa}>0$ and $W,Q\in \mathbb{R}^{n\times n}$ with $W,Q\succ 0$ such that, for all $k\in [M]$, the LMI $\Psi_i(\tau_k)\prec -\varepsilon I_n$ holds. Then, $\Psi_i(t)\prec 0$ for all $t\in[0,T]$.
\end{theorem}
\begin{pf}
For $t\in [0,T]$, by assumption, there exists $j\in [M]$ such that $0<\tau_j-t\leq \mathfrak{s}<\mu$. Employing Corollary~\ref{Cor:PhiDifference}, we conclude that $\left\|\Psi_i(t)-\Psi_i(\tau_j) \right\|<\varepsilon$. By Weyl's inequality \citep[Section 4.3]{horn2012matrix}, the corresponding maximal eigenvalues satisfy 
\begin{equation*}
\left|\lambda_{max}\left(\Psi_i(t) \right)-\lambda_{max}\left(\Psi_i(\tau_j) \right) \right|< \varepsilon.
\end{equation*}
Hence,
$\lambda_{max}\left(\Psi_i(t) \right) <\lambda_{max}\left(\Psi_i(\tau_j) \right) +\varepsilon<0$,
implying that $\Psi_i(t)\prec 0$. Since $t\in [0,T]$ is arbitrary, the proof is concluded. $\hfill\blacksquare$ 
\end{pf}
\begin{remark}
Theorem~\ref{Thm:SamplingLMIs} provides a sufficient condition for guaranteeing the feasibility of the LMI $\Psi_i(t)\prec 0,\ t\in [0,T]$ via discretisation. However, in practice, this method will likely yield a conservative bound on the stencil size $\mathfrak{s}$ in \eqref{eq:Stencil}; tighter estimates on $\mathfrak{s}$ should be attainable via direct computation in concrete examples.
\end{remark}

\section{Numerical examples}\label{sec:sim}
We validate the efficacy of the proposed observer-based sampled-data switching law with dwell-time via two numerical examples, which demonstrate that our approach can stabilise switched systems when all the individual modes have unstable matrices $A_i$, and even when the matrices $A_i$ do not admit a Hurwitz convex combination. The proposed switching strategy is also applied to control a double buck-boost converter, confirming its practical effectiveness in real-world applications to power electronics.

Throughout the examples, given system \eqref{eq:switched_sys}, we choose observer gains \AR{by solving the simplified LMIs in \eqref{eq:LMIDesign1} and setting $L_i = Q^{-1}Y_i,\ i\in [\ell]$}, so that Assumption~\ref{ass:symultaneous_LF} is satisfied.
For a chosen $\Pi \in \mathcal{M}_n$\footnote{The selection of $\Pi \in \mathcal{M}_n$ depends on whether a stabilising convex combination exists: if so, $\Pi$ is chosen such that $\lambda^\top \Pi = 0$, where $\lambda$ is the vector of convex coefficients; otherwise, $\Pi$ is interpreted as the infinitesimal generator of a Markov process where $\Pi_{i,j}$ denotes the expected transition rate from subsystem $i$ to $j$ (see \cite{affine:Russo:2024} and references therein).} and $T,\zeta>0$, we find $X_i$ that satisfy \eqref{eq:LyapMetz}, based on which we then compute $P_i(t)$ by integrating \eqref{eq:P_dot}.
Then, treating $h,\alpha>0$ as tuning parameters, we check the existence of $Q,W\succ 0$ and $\gamma>0$ such that $\Psi_i(t)\prec0$ is satisfied $\forall i\in[\ell]$ and $\forall t \in [0,T]$ by verifying the LMIs on a uniform grid with stencil $\mathfrak{s}=0.01$, where $\mathfrak{s}$ is defined in \eqref{eq:Stencil}. We verify that  
the switching law \eqref{eq:switching_law} ensures global asymptotic stabilisation by simulating the closed-loop system.

\begin{figure}[h]
    \centering
\includegraphics[width=0.85\columnwidth]{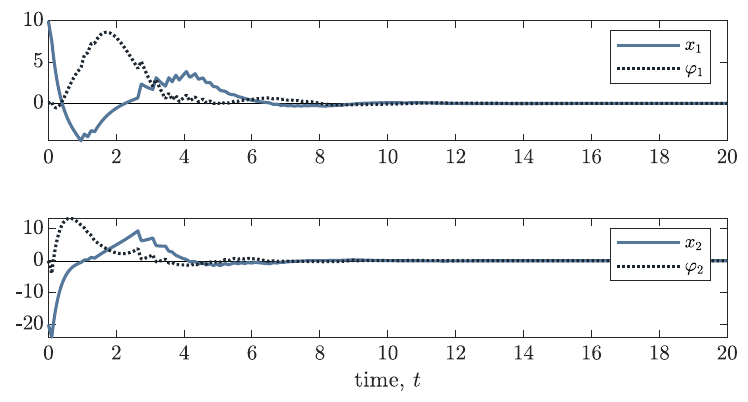}
\vspace{-2mm}
    \caption{\small Time evolution of the state $x$ (blue) and of the observed state $\varphi$ (black) in Example~\ref{sec:example_unstable}.}
\label{fig:unstable_x_phi}
\end{figure}

\begin{figure}[h]
\centering
\includegraphics[width=0.85\columnwidth]{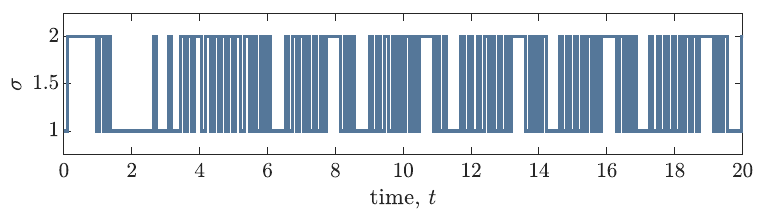}
\vspace{-2mm}
    \caption{\small Time evolution of the switching law $\sigma(t)$ in Example~\ref{sec:example_unstable}.}
\label{fig:unstable_sigma} 
\end{figure}

\begin{example}\label{sec:example_unstable}
\textbf{(Unstable nominal modes.)}
Consider system \eqref{eq:switched_sys} with two modes, $\sigma(t) \in \lbrace 1,\,2\rbrace$, matrices
\begin{align*}
    &A_1 = \begin{bmatrix}
        -2 & 0.3 \\
        -2 & 1
    \end{bmatrix}, \enskip
    A_2 = \begin{bmatrix}
        1 & 2 \\
        -0.3 & -4
    \end{bmatrix}, \\
    & D_1 = D_2 = \begin{bmatrix}
        1 & 1 \\
        1 & -1
    \end{bmatrix},\enskip 
    C = 10^{-5}\cdot \begin{bmatrix}
        1 & 0 \\
        0 & 1
    \end{bmatrix}
\end{align*}
and functions $f_i(x) = \RK{\kappa} \tfrac{\|x\|}{1+\|x\|} \begin{bmatrix}
    1 & 1
\end{bmatrix}^\top$, with $\kappa=0.002$. Both $A_1$ and $A_2$ are unstable. However, their convex combination $A_\lambda = \lambda A_1 + (1-\lambda)A_2$, with $\lambda = 0.5$, is a Hurwitz matrix. We choose observer gains
\begin{equation*}
    L_1 = \begin{bmatrix}
        -1.50 &  -0.85 \\
        -0.85  &  1.50
    \end{bmatrix}, \quad
    L_2 = \begin{bmatrix}
        1.50 &  0.85 \\
        0.85  &  -3.50
    \end{bmatrix}
\end{equation*}
that satisfy Assumption~\ref{ass:symultaneous_LF} and select matrix $\Pi$ in \eqref{eq:LyapMetz} as
\begin{equation*}
    \Pi =  \begin{bmatrix}
        -21.21 & 21.21 \\
        21.21 & -21.21
    \end{bmatrix}.
\end{equation*}
The simulation results for the closed-loop system with the switching law \eqref{eq:switching_law}, inter-sampling interval $h=0.05$, dwell time $T=0.1$, parameters $\zeta = 0.1$ and $\alpha = 10^{-6}$, are shown in Fig.~\ref{fig:unstable_x_phi}. The closed-loop system trajectories converge to the origin, despite the instability of the two nominal modes and the unknown nonlinear term. The switching law \eqref{eq:switching_law}, shown in Fig.~\ref{fig:unstable_sigma}, satisfies the dwell-time constraint.
\end{example}

\begin{figure}[t]
   \centering
\includegraphics[width=0.85\columnwidth]{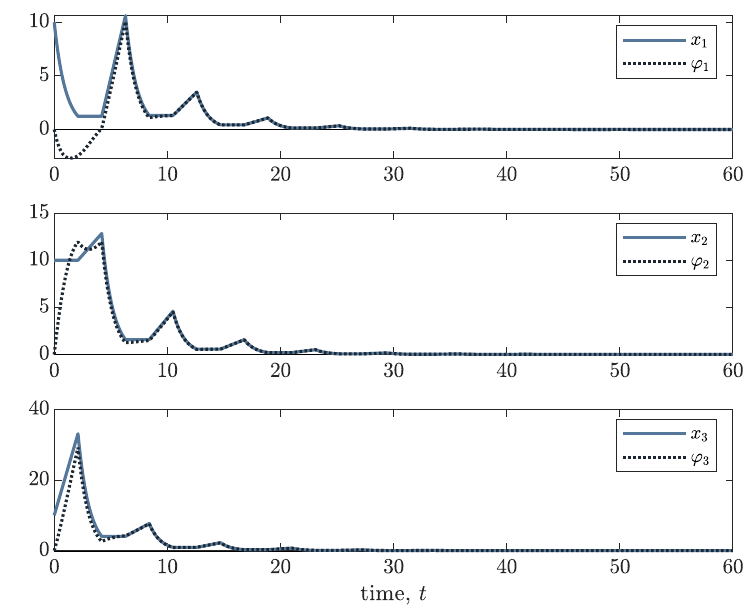}
\vspace{-2mm}    \caption{\small Time evolution of the state $x$ (blue) and of the observed state $\varphi$ (black) in Example~\ref{sec:example_noHurwitz}.}
\label{fig:noHurwitz_x_phi}
\end{figure}

\begin{figure}[t]
    \centering
\includegraphics[width=0.85\columnwidth]{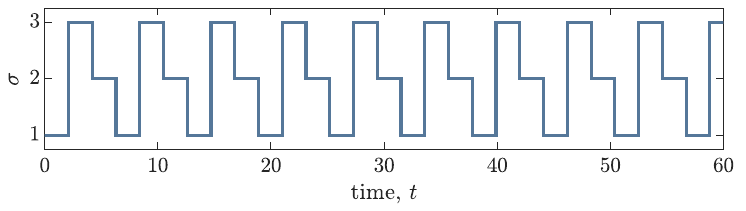}
\vspace{-2mm}    \caption{\small Time evolution of the switching law $\sigma(t)$ in Example~\ref{sec:example_noHurwitz}.}
\label{fig:noHurwitz_sigma}
\end{figure}

\begin{example}\label{sec:example_noHurwitz}
\textbf{(Nominal modes with no Hurwitz convex combination.)}

Consider a switched system as in \eqref{eq:switched_sys} with three modes, $\sigma(t)\in\{1,2,3\}$, matrices 
\begin{align*}
    &A_1 = \begin{bmatrix}
        -\eta & 0 & 0 \\
        0 & 0 & 0 \\
        0 & \beta & 0
    \end{bmatrix}, \enskip
    A_2 = \begin{bmatrix}
        0 & 0 & \beta \\
        0 & -\eta & 0 \\
        0 & 0 & 0
    \end{bmatrix}, \enskip
    A_3 = \begin{bmatrix}
        0 & 0 & 0 \\
        \beta & 0 & 0 \\
        0 & 0 & -\eta
    \end{bmatrix},\\
    & D_1 = D_2 = D_3 = \begin{bmatrix}
        1 & 1 & 1 
    \end{bmatrix},\enskip 
    C =  \begin{bmatrix}
        1 & 0 & 0  \\
        0 & 1 & 0 \\
        0 & 0 & 1
    \end{bmatrix}
\end{align*}
and functions $f_i(x) = \RK{\kappa}\begin{bmatrix}
    \sin(x_1) & \sin(x_2) & \sin(x_3)
\end{bmatrix}^\top$, with $\kappa=0.002$. Matrices $A_i$ are taken from a well-known congestion example (see, e.g., \citealp{traffic:Blanchini:2012}) and are characterised by the absence of a Hurwitz convex combination. Hence, switching approaches such as the ones by \cite{quadratic:bolzern:2004,practical:sanchez:2019,time-triggered:albea:2021} cannot be applied to stabilise the system. We select the system parameters as $\eta=1$ and $\beta = 1.1$ 

We choose observer gains \begin{equation*}
    L_1 = \begin{bmatrix}
        -0.51 \\
        0.53 \\
        0.53 \\
    \end{bmatrix}, \enskip
    L_2 = \begin{bmatrix}
        0.53 \\
        -0.51 \\
        0.53 \\
    \end{bmatrix},
\end{equation*}
that satisfy Assumption~\ref{ass:symultaneous_LF} and select matrix $\Pi$ in \eqref{eq:LyapMetz} as
\begin{equation*}
    \Pi = \begin{bmatrix}
        -10    & 0  &  10 \\
        10 &  -10   &  0\\
        0   & 10 &  -10\\
    \end{bmatrix}.   
\end{equation*}
The simulation results for the closed-loop system with the switching law in Theorem~\ref{thm:main}, inter-sampling interval $h=0.2$, dwell time $T=2.1$,  parameters $\zeta = 0.1$ and $\alpha = 10^{-6}$, are shown in Fig.~\ref{fig:noHurwitz_x_phi}, while the switching law is shown in Fig.~\ref{fig:noHurwitz_sigma}.
The switching law \eqref{eq:switching_law} stabilises the closed-loop system and satisfies the dwell-time constraint.
\end{example}

\begin{example}\label{sec:double_buckboost}
\textbf{(Double buck-boost converter.)}
\begin{figure}
    \centering
	\ctikzset{bipoles/length=4.75cm, font=\fontsize{55}{0}\selectfont}
	\begin{circuitikz} [scale=0.165, transform shape]
    \draw (0,0) to[short] (40,0);
    
    \draw (0,12) to[battery, v_=$\mathcal{E}$, voltage/american label distance=1.8] (0,0);
    \draw (0,12) to[short] (4,12); 
    
    \draw (4,12) to[short, -o] (6,12); 
    \draw (7,9) to[short, -o] (8.5,10.5);  
    \node at (9,8.5) (u1){$u_1$};
    \node at (6,13.5) (u10){$0$};
    \node at (9,13.5) (u11){$1$};
    
    \draw (7,0) to[L, l_=$\mathcal{L}_1$, i=$x_1 \mkern 5mu$] (7,9);
    
    \draw (14,12) to[short, -o] (9,12);
    \draw (14,12) to[C=$\mathcal{C}_1$, v_=$x_2$] (14,0);
    \draw (14,12) to[short] (20,12);
    \draw (20,12) to[R=$\mathcal{R}_1$] (20,0);
    
    \draw (20,12) to[short] (24,12);
    
    \draw (24,12) to[short, -o] (26,12);
    \draw (27,9) to[short, -o] (28.5,10.5);
    \node at (29,8.5) (u2){$u_2$};
    \node at (26,13.5) (u20){$0$};
    \node at (29,13.5) (u21){$1$};
    
    \draw (27,0) to[L, l_=$\mathcal{L}_2$, i=$x_3 \mkern 5mu$] (27,9);
    
    \draw (34,12) to[short, -o] (29,12);
    \draw (34,12) to[C=$\mathcal{C}_2$, v_=$x_4$] (34,0);
    \draw (34,12) to[short] (40,12);
    \draw (40,12) to[R=$\mathcal{R}_L$] (40,0);
\end{circuitikz}
    \caption{Topology of a double buck-boost converter.}
    \label{fig:buck-boost}
\end{figure}

Double buck-boost converters are switching power electronic devices recognised for their ability to independently manage two loads with distinct steady-state or tracking requirements \cite[Section 2.11]{power_conv:SiraRamirez:2006}. As illustrated in Fig.~\ref{fig:buck-boost}, the topology incorporates two switching elements, yielding four distinct operating configurations. Consequently, the converter can be modelled as the switched affine system  \eqref{eq:switched_sys} with $\theta=1$, $f_i\equiv  0$, $\sigma(t) \in [4]$, and
\begin{alignat*}{2}
    & A_1 = \small \begin{bmatrix}
        0 & \mathcal{L}_1^{-1} & 0 & 0 \\
        -\mathcal{C}_1^{-1} & - \left( \mathcal{R}_1\mathcal{C}_1\right)^{-1} & 0 & 0 \\
        0 & 0 & 0 & \mathcal{L}_2^{-1} \\
        0 & 0 & -\mathcal{C}_2^{-1} & -\left(\mathcal{R}_L\mathcal{C}_2 \right)^{-1}
    \end{bmatrix}, \ b_1 = \small \begin{bmatrix}
        0\\0\\0\\0
    \end{bmatrix},\\
    & A_2 = \small \begin{bmatrix}
        0 & 0 & 0 & 0 \\
        0 & -\left(\mathcal{R}_1\mathcal{C}_1 \right)^{-1} & 0 & 0 \\
        0 & 0 & 0 & \mathcal{L}_2^{-1} \\
        0 & 0 & -\mathcal{C}_2^{-1} & -\left( \mathcal{R}_L\mathcal{C}_2\right)^{-1}
    \end{bmatrix}, \ b_2=\begin{bmatrix}
        \mathcal{E}\mathcal{L}_1^{-1} \\
        0 \\
        0 \\
        0
    \end{bmatrix}, \\
    & A_3 = \small \begin{bmatrix}
        0 & \mathcal{L}_1^{-1} & 0 & 0 \\
        -\mathcal{C}_1^{-1} & -\left( \mathcal{R}_1\mathcal{C}_1\right)^{-1} & -\mathcal{C}_1^{-1} & 0 \\
        0 & \mathcal{L}_2^{-1} & 0 & 0 \\
        0 & 0 & 0 & -\left( \mathcal{R}_L\mathcal{C}_2\right)^{-1}
    \end{bmatrix}, \ b_3=b_1,\\
    & A_4 = \small \begin{bmatrix}
        0 & 0 & 0 & 0 \\
        0 & -\left(\mathcal{R}_1\mathcal{C}_1 \right)^{-1} & -\mathcal{C}_1^{-1} & 0 \\
        0 & \mathcal{L}_2^{-1} & 0 & 0 \\
        0 & 0 & 0 & -\left( \mathcal{R}_L\mathcal{C}_2\right)^{-1}
    \end{bmatrix}, \ b_4=b_2,
\end{alignat*}
In this formulation, the state variables are defined as follows: $x_1$ and $x_3$ represent the currents through inductors $\mathcal{L}_1$ and $\mathcal{L}_2$, respectively. The voltage variables $x_2$ and $x_4$ correspond to the capacitor voltages at each conversion stage; specifically, $x_2$ is the voltage across $\mathcal{C}_1$ at the first stage, while $x_4$ represents the voltage across $\mathcal{C}_2$, which serves as the final output or load voltage.
The switching signal $\sigma(t)$ denotes the configurations of the pair of switches $(u_1,u_2)$ in Fig.~\ref{fig:buck-boost}, with $\sigma(t)=1$ corresponding to the configuration  $(0,0)$, $\sigma(t)=2$ to $(1,0)$, $\sigma(t)=3$ to $(0,1)$ and, finally, $\sigma(t)=4$ to $(1,1)$. The objective of the double buck-boost converter is to regulate the voltages at the two output stages, specifically $x_2$ and $x_4$, to their respective reference values, $x_2^\star$ and $x_4^\star$. Hence, for given $x_2^\star$ and $x_4^\star$, a parametrisation of the equilibrium point in terms of the steady-state output voltages is obtained as
\begin{equation*}
    x_1^\star =  \left(\frac{x_2^\star}{\mathcal{R}_1}+\frac{x_4^{\star^2}}{x_2 \mathcal{R}_L}\right)\frac{x_2^\star -\mathcal{E}}{\mathcal{E}}, \quad x_3^\star = \frac{x_4^\star(x_4^\star-x_2^\star)}{\mathcal{R}_L x_2^\star}.
\end{equation*}
Therefore, the control objective is to regulate to zero the error variable $x^e = x - x^\star$,  $x^\star = [x_1^\star\,\,x_2^\star\,\,x_3^\star\,\,x_4^\star]^{\top}$, obeying the switched system
\begin{equation}\label{eq:buckboost_shifted}
    \begin{array}{lll}
        &\dot{x}^e(t) = A_{\sigma(t)}x^e(t)+b^e_{\sigma(t)},\\ &y^e(t) = D_{\sigma(t)} x^e(s_k), \\
        &z^e(t) = Cx^e(t),
    \end{array}
\end{equation}
where $b^e_{\sigma(t)} = A_{\sigma(t)}x^\star+b_{\sigma(t)}$. The numerical values for the converter parameters are selected as in \cite[Section 2.10]{power_conv:SiraRamirez:2006}. Given the input voltage $\mathcal{E}=12V$, with $x_2^\star = -12V$ and $x_4^\star = 12V$, the convex combination of $\left\{A_i \right\}_{i=1}^4$ obtained with coefficients $\lambda_1= \lambda_4 = 0.5$, $\lambda_2 = \lambda_3 = 0$, results in a stable matrix. We consider our system with $D_i=I_4,\ i\in [4]$, which allows for comparison with \cite{sampled:Hetel:2013}, and we choose the observer gains
\begin{equation*}
    \begin{array}{c}
    L_1 = 10^{-6}\cdot \small \begin{bmatrix}
    0.500 &   -0.042 &   0.002  & -0.019 \\
    -0.042 & 160.105 &  -0.038  &  -0.016 \\
    0.002 &  -0.038  &  0.204  & -0.301 \\
    -0.019 &  -0.016 &   -0.301  & 64.623
    \end{bmatrix},\vspace{0.1cm} \\
    L_2 = 10^{-6}\cdot \small \begin{bmatrix}
    0.750 &  -0.064  &  0.003 &  -0.028 \\
    -0.064 &  240.158 &  -0.058 &  -0.024 \\
    0.003  & -0.058 &   0.307 &   -0.452 \\
    -0.028  & -0.024  & -0.452  & 96.935
    \end{bmatrix},\vspace{0.1cm}  \\
    L_3 = 10^{-6}\cdot \small \begin{bmatrix}
    1.000  & -0.085 &   0.004 &  -0.038 \\
    -0.085 & 320.211 &  -0.077 &   -0.032 \\
    0.004 &  -0.077  &  0.409 &  -0.603 \\
    -0.038  & -0.032  &  -0.603 &  129.247
    \end{bmatrix},  \vspace{0.1cm}  \\
    L_4 = 10^{-6}\cdot\small \begin{bmatrix}
    1.250 &  -0.106 &    0.005 &  -0.047 \\
    -0.106 &  400.263 &  -0.096  & -0.040 \\
    0.005  & -0.096 &    0.511 &  -0.753 \\
    -0.047 &  -0.040  & -0.753 & 161.558
    \end{bmatrix},
    \end{array}
\end{equation*}
that satisfy Assumption~\ref{ass:symultaneous_LF}, as well as the matrix $\Pi$ in \eqref{eq:LyapMetz} as
\begin{equation*}
    \Pi = 10^8\cdot \small \begin{bmatrix}
    -0.7071    &    0   &      0 &   0.7071 \\
    1  & -1  &       0  &       0 \\
    1     &    0  &  -1    &     0 \\
    0.7071     &    0    &     0 &   -0.7071
    \end{bmatrix}.
\end{equation*}
Finally, we set $\zeta = 10^{-3}$ and $\alpha = 10^{-4}$.
Due to the high switching frequency required for DC-DC power converters, both the dwell time and the maximum sampling interval are set to $T=10^{-5}$ s  and $h=10^{-5}$ s, while the converter is assumed to be initially discharged, i.e., $x(0) = 0$.  In this example, it is assumed that the affine vectors $b_i$ are known, and it is hence possible to implement the switching strategy outlined in Section~\ref{Sec:ModifCLaw} and summarised in Theorem~\ref{thm:mainMod}. The proposed switching law is compared with the algorithm proposed by \cite{sampled:Hetel:2013}, which addresses the state feedback stabilisation of switched affine systems under sampled-data control. To ensure a fair comparison between the two approaches, the sampling interval for the algorithm by \cite{sampled:Hetel:2013} is set to $h$, matching the sampling interval length used in our proposed approach. The performance of our proposed switching law is shown in Fig.~\ref{fig:double_buckboost_x}, demonstrating that the states converge towards the desired reference values and oscillate around them. The obtained performance is comparable to the performance obtained from the approach by \cite{sampled:Hetel:2013}, in terms of both transient and steady-state dynamics. \AR{In fact, for the two capacitor voltages $x_2$ and $x_4$, our switching strategy maintains a steady-state relative error below $0.76\%$ and $0.18\%$, respectively. In comparison, the method proposed by \cite{sampled:Hetel:2013} yields steady-state relative errors below $0.61\%$ and $0.10\%$ for $x_2$ and $x_4$, respectively. Furthermore, the proposed algorithm yields a steady-state RMS error equal to  $0.047$~V and $0.010$~V for $x_2$ and $x_4$, respectively, while the algorithm presented in \cite{sampled:Hetel:2013} guarantees a steady-state RMS error equal to $0.042$~V and $0.006$~V, for $x_2$ and $x_4$, respectively.} As shown in Fig.~\ref{fig:double_buckboost_sigma}, the switching signals obtained by both algorithms eventually converge to a repeated switching between modes 1 and 4. 

While the design of the switching strategy by \cite{sampled:Hetel:2013} imposes a dwell-time constraint that is equal to the sampling period, our proposed methodology decouples the dwell time from the sampling time. This separation allows multiple switching events within a single sampling interval, a feature particularly useful when the sampling and switching are not collocated and the measurements are acquired at a lower frequency than the maximum achievable switching frequency. This advantage is clearly demonstrated in Fig.~\ref{fig:double_buckboost_smallerDT}, where both algorithms are implemented with a sampling period of $10^{-5}s$ , while our approach employs a shorter dwell time of $10^{-6}s$ . As illustrated in the figure, our proposed algorithm achieves significantly lower voltage oscillation amplitudes around the desired set-points, due to the refined dwell-time control. \AR{In fact, in this case our approach yields a steady-state relative error below $0.067\%$ and $0.012\%$ for the two capacitor voltages, respectively. Furthermore, adopting smaller dwell time guarantees a steady-state RMS error equal to $0.0046$~V and $0.0007$~V, for $x_2$ and $x_4$, respectively.}

\begin{figure}
    \centering
    \includegraphics[width=\columnwidth]{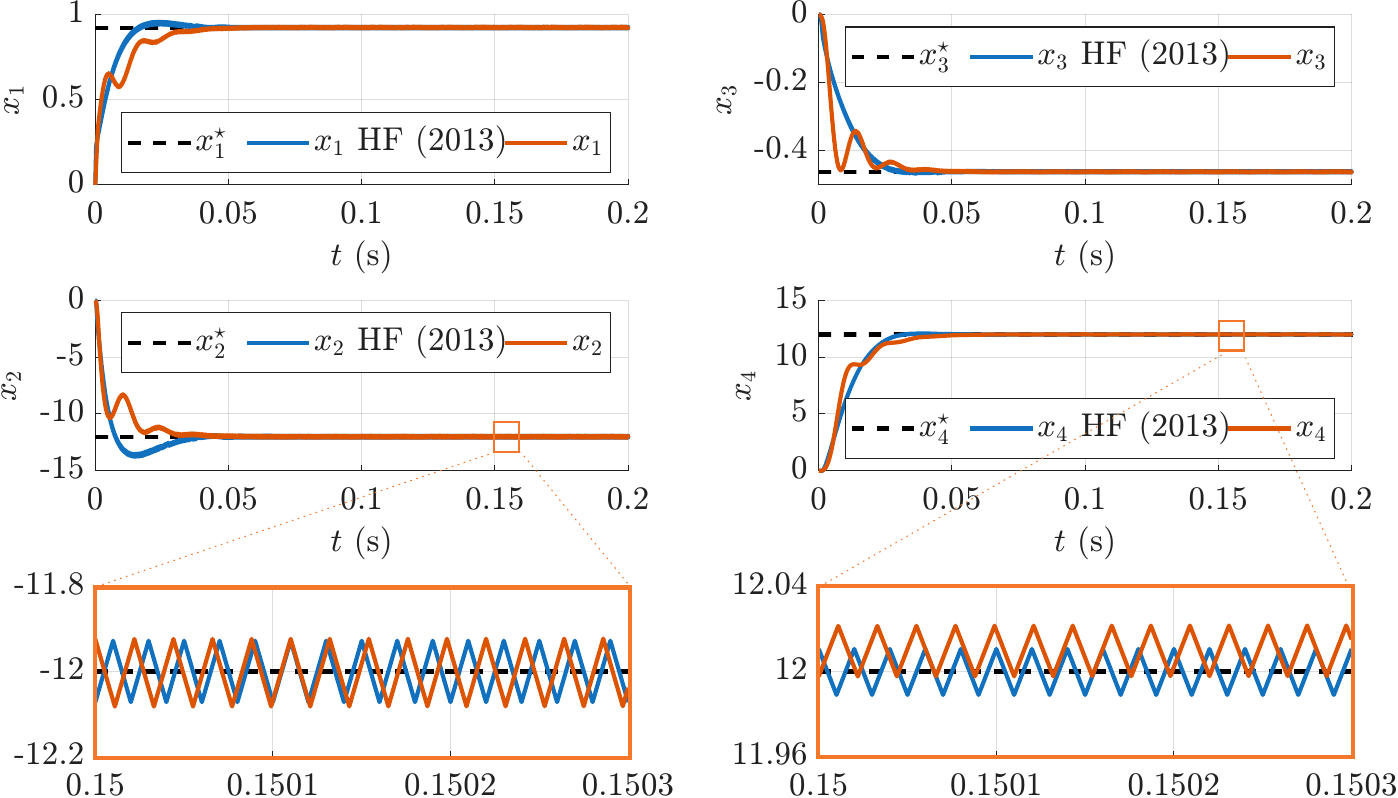}
    \caption{Time-evolution of the double buck-boost converter states obtained with our proposed switching strategy (red solid line) and the strategy by \cite{sampled:Hetel:2013} (blue solid line); the black dashed line is the reference signal. Top left: current $x_1(t)$. Bottom left: first stage voltage $x_2(t)$. Top right: current $x_3(t)$. Bottom right: second stage voltage $x_4(t)$.}
    \label{fig:double_buckboost_x}
\end{figure}

\begin{figure}
    \centering
    \includegraphics[width=0.85\columnwidth]{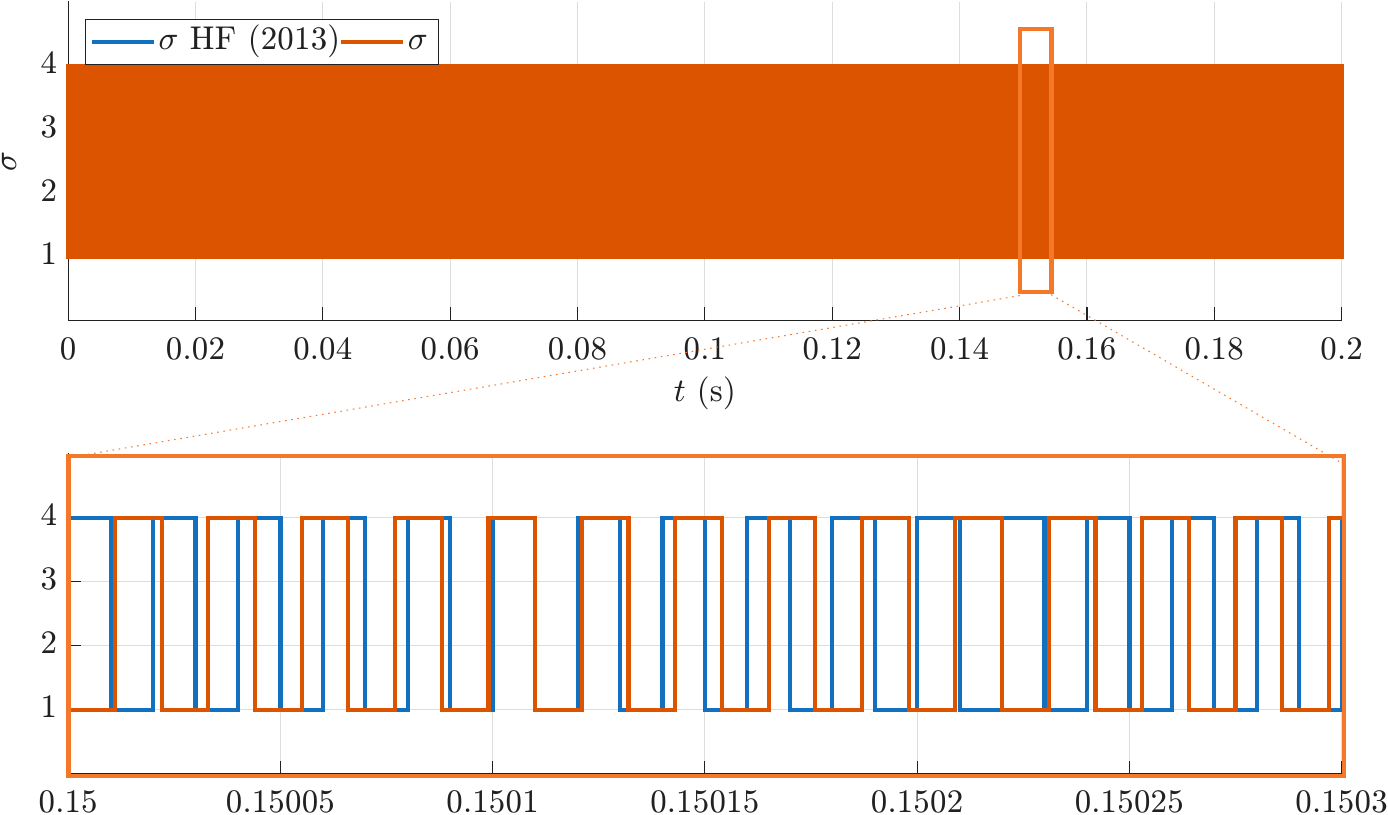}
    \caption{Time-evolution of the double buck-boost converter switching signal obtained with the proposed switching strategy (red solid line) and the strategy in \cite{sampled:Hetel:2013} (blue solid line).}
    \label{fig:double_buckboost_sigma}
\end{figure}

\begin{figure}
    \centering
    \includegraphics[width=\columnwidth]{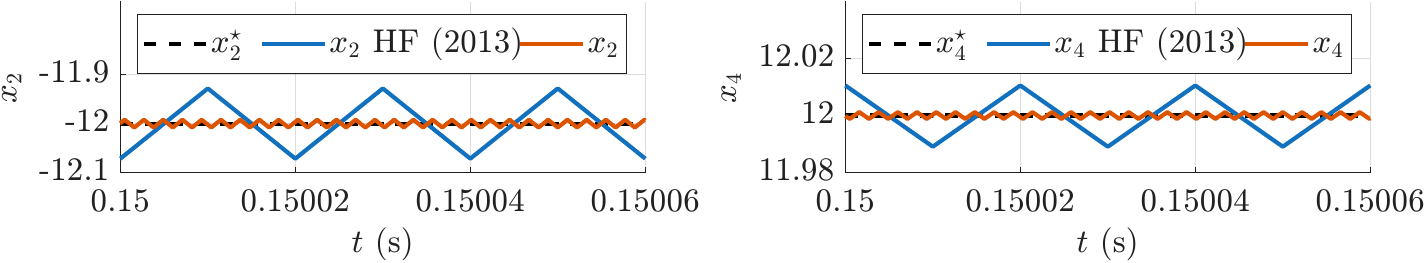}
    \caption{Ultimate bound of the double buck-boost converter voltages obtained with our proposed switching strategy (red solid line) and the strategy by \cite{sampled:Hetel:2013} (blue solid line); the black dashed line is the reference signal. In our approach, the dwell time is smaller than the sampling time, whereas the strategy by \cite{sampled:Hetel:2013} forces the sampling period to be equal to the dwell-time constraint.}
    \label{fig:double_buckboost_smallerDT}
\end{figure}

\end{example}

\section{Concluding Discussion}\label{sec:conclusions}
We proposed a novel observer-based switching law with dwell-time constraints to achieve either stability or ultimate boundedness of switched systems with nominal linear-affine dynamics, subject to uncertain Lipschitz nonlinearities and sampled-data output measurements. Our switching law relies on a suitably designed observer and Lyapunov-Metzler inequalities. We performed Lyapunov stability analysis, leading to time-dependent LMI conditions that provide a lower bound on the maximum inter-sampling time, and of the Lipschitz constants $\RK{\kappa}$, for which the closed-loop system is globally asymptotically stable (or, in the presence of switching affine terms, ultimately bounded) with the designed switching law. 
We investigated the feasibility of the derived LMIs and (i) showed how they can be modified to yield observer gains design; (ii) provided equivalent reduced-order LMI conditions; and (iii) proved that the time dependence of the LMIs can be removed by discretisation on a finite grid. Numerical simulations confirmed our theoretical results. We also compared our method to an existing approach within the literature, by \cite{sampled:Hetel:2013}, which deals with sampled-data state-feedback stabilisation, without including observers and without allowing for dwell-time constraints independent of the inter-sampling time of the state.

Although, to the best of our knowledge, no other work addresses the same framework of switched affine systems with nonlinearities under dwell-time constraints, sampled-data measurements and output feedback, other contributions in the literature share similarities with specific aspects of our problem formulation, but with important differences. For instance, \cite{stabilization:Perez:2015} consider the robust stabilisation of nonlinear switched systems affected by bounded perturbations and by a control input $u(t)$, and propose an observer-based switching controller using sampled and quantised output measurements. As in our case, the nonlinearities satisfy a quasi-Lipschitz condition. A dwell-time constraint is also considered, and sufficient conditions for practical stability are formulated through Bilinear Matrix Inequalities derived from an extended invariant ellipsoid method. However, the approach by \cite{stabilization:Perez:2015} differs from ours in several key aspects: their switching law is time-dependent, with in-built avoidance of Zeno behaviour via the introduction of a dwell-time, whereas we design an observer-dependent switching signal \eqref{eq:switching_law}; their practical stability result relies on the presence of a control input $u(t)$ in the dynamics (which relies on the designed observer), whereas in our case stabilisation/ultimate boundedness is achieved through the observer-based switching signal design only; their stability analysis employs a Lyapunov functional suitable for general time-delay systems, more conservative than the Wirtinger-based functional we use; they do not employ Lyapunov-Metzler inequalities in their design.\\
\cite{stabilization:Li:2021} consider interval observer-based controller design for switched linear (but not linear-affine) systems, with Lipschitz continuous disturbance and measurement noise, relying on the observer designed by \citep{ethabet2018interval}, but neither of these works considers sampled-data measurements. Also, in both works, the switching signal is time-dependent, with the avoidance of Zeno behaviour imposed via an appropriate dwell time. \cite{stabilization:Li:2021} introduce a controller $u(t)$ into the dynamics to achieve stability under \emph{average} dwell time constraints. Moreover, both works impose structural assumptions on the system measurements by requiring that there exist observer gains $\left\{L_i \right\}_{i\in [\ell]}$ such that $A_i-L_iC_i$ is Metzler for all $i\in [\ell]$, in order to apply techniques from cooperative systems.\\
\cite{stabilization:Zhao:2020} propose an observer-based reliable control framework for switched linear (but not linear-affine) systems; the measurements are performed in continuous time and the time-dependent switching mechanism is exogenous. An output-feedback event-triggered controller $u(t)$ is designed in order to achieve stability of the closed-loop system; sampling only appears due to the use of an event-triggering mechanism in the controller. The analysis relies on the separation principle for the observer and state dynamics, and cannot be applied in the presence of Lipschitz nonlinearities.\\
Analogously, \cite{stabilization:Fu:2021} consider in Chapter 3 linear (but not linear-affine) switched systems subject to continuous time measurements and an exogenously given time-dependent switching signal, and investigate the stabilisation of the system via an event-triggered observer-based controller $u(t)$; also here, sampling only appears due to the use of an event-triggering mechanism in the controller, and the absence of nonlinearities in the system again allows for separation of observer and controller designs.\\
\cite{stabilization:Zhang:2025} consider only observation of a controlled switched system, subject to exogenously given time-dependent switching signal, sampled-data state measurements and uncertain nonlinearities, with the control objective being the achievement of exponentially ultimately bounded estimating error. The unknown nonlinearity is subject to global Lipschitz conditions, with the Lipschitz bound defining the ultimate bond on the estimation error. The observer is constructed in Luenberger form and ultimate boundedness is shown via Lyapunov analysis leading to sum of squares and LMI conditions. However, they consider neither the stabilisation/ultimate boundedness of the original state equation nor the design of a state/output-based switching signal, and do not employ Lyapunov-Metzler inequalities.\\
Finally, \cite{stabilization:Carneiro:2024} investigate practical stabilisation of a continuous-time uncertain switched affine system, with a switching affine term, subject to continuous time measurements, corrupted by noise. They employ a Luenberger observer whose state is then sampled at discrete time instants and used for synthesis of the switching signal. The switching law design and the stability analysis in the presence of sampled-data by \cite{stabilization:Carneiro:2024} is similar to that by \cite{sampled:Hetel:2013} (in both, the sampling instants coincide with the switching instants, differently from our framework where the dwell time and the sampling time are decoupled) and assumes the existence of a convex combination $A(\kappa)$ of the mode matrices $\left\{A_i \right\}$ that is Hurwitz, as well as the fact that the pairs $\left(A(\kappa),C_i \right)$, where $C_i$ are the measurement operators, are all observable. In comparison, our approach, which also handles the presence of affine switching dependent terms in the dynamics, in addition to Lipschitz nonlinearities, does not require a Hurwitz convex combination $A(\kappa)$ to exist, as also demonstrated in Example~\ref{sec:example_noHurwitz}, where the nominal modes do not admit a Hurwitz convex combination.

Future work includes further extending our approach to consider general time-varying delays and more general classes of switched systems.

\bibliography{ifacconf}

\end{document}